\documentclass{amsart}
\usepackage{amssymb}
\usepackage{amsmath}
\usepackage{amsfonts}

\input{tcilatex}

\begin{document}
\title[Dissimilarity Semimetrics for Preferences]{{\Large \textbf{A Class of
Dissimilarity Semimetrics for Preference Relations}}}
\author{Hiroki Nishimura}
\address{Department of Economics,\\
University of California at Riverside}
\email{hiroki@ucr.edu}
\author{Efe A. Ok}
\address{Department of Economics and Courant Institute of Mathematical
Sciences, New York University}
\email{efe.ok@nyu.edu}
\date{March 08, 2022}
\thanks{This paper is in final form and no version of it will be submitted
for publication elsewhere.}

\begin{abstract}
We propose a class of semimetrics for preference relations any one of which
is an alternative to the classical Kemeny-Snell-Bogart metric. (We take a
fairly general viewpoint about what constitutes a preference relation,
allowing for any acyclic order to act as one.) These semimetrics are based
solely on the implications of preferences for choice behavior, and thus
appear more suitable in economic contexts and choice experiments. In our
main result, we obtain a fairly simple axiomatic characterization for the
class we propose. The apparently most important member of this class (at
least in the case of finite alternative spaces), which we dub \textit{the
top-difference semimetric}, is characterized separately. We also obtain
alternative formulae for it, and relative to this metric, compute the
diameter of the space of complete preferences, as well as the best
transitive extension of a given acyclic preference relation. Finally, we
prove that our \textit{preference} \textit{metric spaces} cannot be
isometically embedded in a Euclidean space.
\end{abstract}

\maketitle

\section{Introduction}

The matter of distinguishing between individual preference relations on a
set of choice objects is of great import for a variety of subdisciplines of
economics, sociology, political science, and psychology. It is often the
case that researchers wish to understand how dissimilar are the preferences
of subjects that are estimated in a choice experiment, thereby getting a
sense of the variability of preferences in the aggregate. Or, depending on
the context, one may wish to have a way of determining which of two
individuals is more altruistic (or resp., patient, or risk averse) by
comparing their preferences to a benchmark altruistic (resp., fully patient,
or risk neutral) preference relation. Similarly, we may try to understand
which of two preference relations exhibits more indecisiveness among
alternatives by checking how far off they are from being a complete
preference relation. Or one may wish to investigate the extent to which a
given preference relation violates a rationality axiom by checking how
distant this relation is from the class of all preferences which satisfy
that axiom.

These types of considerations provide motivation for developing general
methods of making dissimilarity comparisons between the family of all
preference relations on a given finite set $X$ of alternatives.\footnote{%
We are being deliberately loose in this section about what we mean by a
\textquotedblleft preference relation\textquotedblright\ on $X.$ Economists
often take this to mean a preorder (if they wish to allow for
indecisiveness), or a total preorder (if they want to model the preferences
of a decisive individual). By contrast, in voting theory, and operations
research at large, one often assumes indifferences away, and refer to any
partial or linear order as a preference relation. In this paper we work with
acyclic orders, and include all of these specifications as special cases;
see Section 2.2.} The most common way of doing this is by means of equipping
this family with a suitable distance function. The starting point of the
related literature is the seminal work of Kemeny and Snell \cite{KS} who
axiomatically proposed a distance function over linear orders on $X$ -- the
order-theoretic terminology we use in this paper is outlined in Section 2.1
-- which is based on counting the number of rank reversals between two such
orders. (The distance between two linear orders according to this metric is
twice the total number of involved rank reversals.) While its restriction to
linear orders is limiting, Bogart \cite{Bogart1} has extended this metric to
the context of all partial orders on $X$ by means of a modified system of
axioms. To be precise, let us denote the indicator function of any partial
order $\succsim $ on $X$ by $I_{\succsim }$ (that is, $I_{\succsim }$ is the
map on $X\times X$ with $I_{\succsim }(x,y):=1$ if $x\succsim y,$ and $%
I_{\succsim }(x,y):=0$ otherwise). Then, the \textit{Kemeny-Snell-Bogart} 
\textit{metric} on the set of all partial orders on $X$ is defined by%
\begin{equation*}
d_{\text{KSB}}(\succsim ,\succsim ^{\prime })=\sum_{x,y\in X}\left\vert
I_{\succsim }(x,y)-I_{\succsim ^{\prime }}(x,y)\right\vert \text{.}
\end{equation*}%
This metric has been found of great use in deducing a consensus ranking from
a given collection of individual preferences (which may or may not leave
some alternatives unranked). Moreover, the literature provides several
extensions of, and alternatives to, this distance function. (See \cite{Cook2}
for a survey of this literature.)

There are several perspectives in which two preference relations may differ
from each other, and it is of course unreasonable to expect a single
distance function to be sensitive to all of these. Indeed, there is an
aspect, which is of utmost importance for economic analysis, that is not
correctly attended by the Kemeny-Snell-Bogart metric. In economics at large,
a preference relation $\succsim $ is viewed mainly as a means toward making
choices in the context of various menus (i.e., nonempty subsets of the grand
set $X$ with at least two members), where a \textquotedblleft
choice\textquotedblright\ in a menu $S$ on the basis of $\succsim $ is
defined as a maximal element of $S$ with respect to $\succsim $.
Consequently, the more distinct the induced \textquotedblleft
choices\textquotedblright\ of two preference relations across menus are,
there is reason to think of those preferences as being less similar. Here
are two simple examples that highlight in what sense the $d_{\text{KSB}}$
metric does not reflect this viewpoint properly.

\ifx\JPicScale\undefined

\fi

\unitlength.7 mm 
\begin{picture}(65,68)(25,0)

\linethickness{0.2mm}
\put(80,29){\line(0,1){28}}
\put(80,20){\makebox(0,0)[cc]{{\footnotesize {$\succsim$}}}}
\put(80,57){\circle*{2}}
\put(80,50.25){\circle*{2}}
\put(80,43.5){\circle*{2}}
\put(80,36.25){\circle*{2}}
\put(80,29){\circle*{2}}

\put(86,57){\makebox(0,0)[cc]{{\footnotesize {$x_1$}}}}
\put(86,49.5){\makebox(0,0)[cc]{{\footnotesize {$x_2$}}}}
\put(86,43){\makebox(0,0)[cc]{{\footnotesize {$x_3$}}}}
\put(86,36){\makebox(0,0)[cc]{{\footnotesize {$x_4$}}}}
\put(86,28){\makebox(0,0)[cc]{{\footnotesize {$x_5$}}}}

\linethickness{0.2mm}
\put(102,29){\line(0,1){28}}
\put(103,20){\makebox(0,0)[cc]{{\footnotesize {$\succsim_1$}}}}

\put(102,57){\circle*{2}}
\put(102,50.25){\circle*{2}}
\put(102,43.5){\circle*{2}}
\put(102,36.25){\circle*{2}}
\put(102,29){\circle*{2}}

\put(108,57){\makebox(0,0)[cc]{{\footnotesize {$x_2$}}}}
\put(108,49.5){\makebox(0,0)[cc]{{\footnotesize {$x_1$}}}}
\put(108,43){\makebox(0,0)[cc]{{\footnotesize {$x_3$}}}}
\put(108,36){\makebox(0,0)[cc]{{\footnotesize {$x_4$}}}}
\put(108,28){\makebox(0,0)[cc]{{\footnotesize {$x_5$}}}}

\linethickness{0.2mm}
\put(125,29){\line(0,1){28}}
\put(127,20){\makebox(0,0)[cc]{{\footnotesize {$\succsim_2$}}}}

\put(125,57){\circle*{2}}
\put(125,50.25){\circle*{2}}
\put(125,43.5){\circle*{2}}
\put(125,36.25){\circle*{2}}
\put(125,29){\circle*{2}}

\put(131,57){\makebox(0,0)[cc]{{\footnotesize {$x_1$}}}}
\put(131,49.5){\makebox(0,0)[cc]{{\footnotesize {$x_2$}}}}
\put(131,43){\makebox(0,0)[cc]{{\footnotesize {$x_3$}}}}
\put(131,36){\makebox(0,0)[cc]{{\footnotesize {$x_5$}}}}
\put(131,28){\makebox(0,0)[cc]{{\footnotesize {$x_4$}}}}

\put(104,7){\makebox(0,0)[cc]{{\bf {\footnotesize Figure 1}}}}

\end{picture}

\noindent \textit{Example 1.1}. Let $X:=\{x_{1},...,x_{5}\},$ and consider
the linear orders $\succsim ,$ $\succsim _{1}$ and $\succsim _{2}$ on $X$
whose Hasse diagrams are depicted in Figure 1. Clearly, both $\succsim _{1}$
and $\succsim _{2}$ are obtained from $\succsim $ by reversing the ranks of
two alternatives, namely, those of $x_{1}$ and $x_{2}$ in the case of $%
\succsim _{1}$ and those of $x_{4}$ and $x_{5}$ in the case of $\succsim _{2}
$. Consequently, the Kemeny-Snell-Bogart metric judges the distance between $%
\succsim $ and $\succsim _{1}$ and that between $\succsim $ and $\succsim
_{2}$ the same: $d_{\text{KSB}}(\succsim ,\succsim _{1})=2=d_{\text{KSB}%
}(\succsim ,\succsim _{2})$. But this conclusion is not supported from a
choice-theoretic standpoint. Consider an individual whose preferences are
represented by $\succsim $. This person would never choose either $x_{4}$ or 
$x_{5}$ in any menu $S\subseteq X$ with the exception of $S=\{x_{4},x_{5}\}$%
. Consequently, the choice behavior of this person would differ from that of
an individual with preferences $\succsim _{2}$ in only \textit{one} menu,
namely, $\{x_{4},x_{5}\}.$ By contrast, the choice behavior entailed by $%
\succsim $ and $\succsim _{1}$ are distinct in every menu that contains $%
x_{1}$ and $x_{2}$. So if we observed the choices made by two people with
preferences $\succsim $ and $\succsim _{1},$ we would see them make
different choices in \textit{eight} separate menus. From the perspective of
induced choice behavior, then, it is only natural that we classify
\textquotedblleft $\succsim $ and $\succsim _{1}$\textquotedblright\ as
being less similar than \textquotedblleft $\succsim $ and $\succsim _{2}$%
.\textquotedblright \footnote{%
This viewpoint is also advanced in a few other papers in the literature,
namely, Can \cite{Can}, Hassanzadeh and Milenkovic \cite{Has-M}, and Klamler 
\cite{Klamler}. We will clarify the connections between these papers and the
present one as we proceed.}'\footnote{%
As we mentioned above, there are some well-known alternatives to $d_{\text{%
KSB}}$, such as the metrics of Blin \cite{Blin} and Cook and Seiford \cite%
{Cook-S}$.$ These variants are also based on the idea of counting the rank
reversals between two preferences in one way or another, and also yield the
same conclusion as $d_{\text{KSB}}$ in the context of this example.} $\Vert $

\bigskip

This example points to the fact that, at least from the perspective of
choice behavior, the dissimilarity of two preferences depends not only on
the number of rank reversals between them, but also \textit{where} those
reversals occur.\footnote{%
To put this point in a concrete perspective, recall that in the 2020 U.S.
presidental elections, there were four candidates in the Electoral College:
(1) D. Trump and M. Pence, (2) J. Biden and K. Harris, (3) H. Hawkins and A.
N. Walker, (4) J. Jorgensen and S. Cohen. Now consider four voters each
putting candidates (1) and (2) above the candidates (3) and (4). Suppose two
of these voters disagree between the ranking of Trump-Pence and
Biden-Harris, but agree on the relative ranking of (3) and (4), while the
other two are both Trump supporters who happen to disagree on the relative
ranking of (3) and (4). Obviously, in the elections, the latter two
individuals both voted for the Trump-Pence ticket, while the former two
casted opposite votes. And yet the Kemeny-Snell-Bogart metric views the
preferences of these two pairs of voters equally distant from each other!}

In the next example, we illustrate that the Kemeny-Snell-Bogart metric
behaves in a counterintuitive fashion (from the standpoint of induced choice
behavior) also when we allow for non-comparability, or indifference, of some
alternatives.

\bigskip

\noindent \textit{Example 1.2}. Let $X:=\{x_{1},...,x_{4}\},$ and consider
the partial orders $\succsim ,$ $\succsim _{1}$ and $\succsim _{2}$ on $X$
whose Hasse diagrams are depicted in Figure 2. Here $\succsim _{1}$ is
obtained from $\succsim $ by reversing the ranks of the second-best and
worst alternatives, namely, those of $x_{2}$ and $x_{4}$; we have $d_{\text{%
KSB}}(\succsim ,\succsim _{1})=6.$ On the other hand, the third preference $%
\succsim _{2}$ seems very different than $\succsim $ in that it cannot
render a judgement about the relative desirability of \textit{any}
alternative; this is the preference relation of a person who is entirely
indecisive about the alternatives $x_{1},...,x_{4}$ (whatever may be their
reasons). And yet we again have $d_{\text{KSB}}(\succsim ,\succsim _{2})=6.$
This is, again, difficult to accept from a choice-theoretic perspective. The
choices made on the bases of $\succsim $ and $\succsim _{1}$ differ from
each other in exactly four menus. By contrast, there is no telling as to the
precise nature of choices on the basis of $\succsim _{2}$ as every
alternative in every menu is maximal with respect to this relation, so we
have to declare all alternatives on a menu as a potential choice relative to
this preference relation. But then, the choices induced by $\succsim $ and $%
\succsim _{2}$ differ at every menu.\footnote{%
A similar conclusion would hold if the third preference here declared all
alternatives indifferent (instead of incomparable). In that case, the
standard modification of $d_{\text{KSB}}$ would be defined the same way but
with $I_{\succsim }(x,y):=1$ if $x\succ y$ and $I_{\succsim }(x,y):=1/2$ if $%
x\sim y$ (where $\succ $ and $\sim $ are the asymmetric and symmetric parts
of $\succsim ,$ respectively), and this modified metric would judge $%
\succsim $ and $\succsim _{1},$ and $\succsim $ and $\succsim _{2}$ (where
now $x_{1}\sim _{2}x_{2}\sim _{2}x_{3}\sim _{2}x_{4}$) equally distant, even
though the choices induced by $\succsim $ and $\succsim _{2}$ are distinct
from each other at every menu.} $\Vert $

\ifx\JPicScale\undefined

\fi

\unitlength.7 mm 
\begin{picture}(65,68)(25,0)

\linethickness{0.2mm}
\put(80,29){\line(0,1){28}}
\put(80,20){\makebox(0,0)[cc]{{\footnotesize {$\succsim$}}}}
\put(80,57){\circle*{2}}
\put(80,47){\circle*{2}}
\put(80,38){\circle*{2}}
\put(80,29){\circle*{2}}

\put(86,57){\makebox(0,0)[cc]{{\footnotesize {$x_1$}}}}
\put(86,47){\makebox(0,0)[cc]{{\footnotesize {$x_2$}}}}
\put(86,38){\makebox(0,0)[cc]{{\footnotesize {$x_3$}}}}
\put(86,28){\makebox(0,0)[cc]{{\footnotesize {$x_4$}}}}

\linethickness{0.2mm}
\put(102,29){\line(0,1){28}}
\put(103,20){\makebox(0,0)[cc]{{\footnotesize {$\succsim_1$}}}}

\put(102,57){\circle*{2}}
\put(102,47){\circle*{2}}
\put(102,38){\circle*{2}}
\put(102,29){\circle*{2}}

\put(108,57){\makebox(0,0)[cc]{{\footnotesize {$x_1$}}}}
\put(108,47){\makebox(0,0)[cc]{{\footnotesize {$x_4$}}}}
\put(108,38){\makebox(0,0)[cc]{{\footnotesize {$x_3$}}}}
\put(108,28){\makebox(0,0)[cc]{{\footnotesize {$x_2$}}}}

\linethickness{0.2mm}
\put(136,20){\makebox(0,0)[cc]{{\footnotesize {$\succsim_2$}}}}

\put(125,43.5){\circle*{2}}
\put(132,43.5){\circle*{2}}
\put(139,43.5){\circle*{2}}
\put(146,43.5){\circle*{2}}

\put(125,38){\makebox(0,0)[cc]{{\footnotesize {$x_1$}}}}
\put(132,38){\makebox(0,0)[cc]{{\footnotesize {$x_2$}}}}
\put(139,38){\makebox(0,0)[cc]{{\footnotesize {$x_3$}}}}
\put(146,38){\makebox(0,0)[cc]{{\footnotesize {$x_4$}}}}

\put(104,7){\makebox(0,0)[cc]{{\bf {\footnotesize Figure 2}}}}

\end{picture}

These examples demonstrate that there is room for looking at alternatives to
the Kemeny-Snell-Bogart metric and its variants, especially if we wish to
distinguish between preferences on the basis of their implications for
choice. Our proposal here is to define a class of such alternatives looking
directly at the size of the differences in choices induced by preferences
across all menus, where by a \textquotedblleft choice induced by a
preference in a menu $S$,\textquotedblright\ we mean, as usual, any maximal
element in $S$ relative to that preference. So, on a given menu $S,$ we
propose to capture the dissimilarity of two preference relations on $X,$
say, $\succsim $ and $\trianglerighteq $, by comparing the set of $%
M(S,\succsim )$ of all $\succsim $-maximal elements in $S$ with the set $%
M(S,\trianglerighteq )$ of all $\trianglerighteq $-maximal elements in $S$.
A particularly simple way of making this comparison is, of course, just by
counting the elements in $M(S,\succsim )$ that are not in $%
M(S,\trianglerighteq ),$ as well as those in $M(S,\trianglerighteq )$ that
are not in $M(S,\succsim )$. Thus, the number of elements in the symmetric
difference $M(S,\succsim )\triangle M(S,\trianglerighteq )$ tells us how
different $\succsim $ and $\trianglerighteq $ are in terms of the choice
behavior they entail at the menu $S.$ Then, summing over all menus yields
the main semimetric $D$ we propose here:%
\begin{equation*}
D(\succsim ,\trianglerighteq )=\sum_{S\subseteq X}\left\vert M(S,\succsim
)\triangle M(S,\trianglerighteq )\right\vert \text{.}
\end{equation*}%
We call this map the \textit{top-difference semimetric. }

A reinterpretion of this semimetric by using choice theory is in order. Let
us first recall that a \textit{choice correspondence} on $X$ is any function 
$C:2^{X}\rightarrow 2^{X}$ with $C(S)\subseteq S.$ If we abtract away from
how choice correspondences come to being (via preference maximization, or
boundedly rational choice procedures, or randomizations, etc.), and treat
them as set-valued functions on the finite set $2^{X}$, then the natural $%
\ell _{1}$-type metric on the set of all choice correspondences on $X$ is of
the form%
\begin{equation*}
d_{\text{K}}(C,C^{\prime })=\sum_{S\subseteq X}\left\vert C(S)\triangle
C^{\prime }(S)\right\vert \text{.}
\end{equation*}%
This metric was indeed proposed, and axiomatically characterized, by Klamler 
\cite{Klamler} (which is why we denote it by $d_{\text{K}}$). Now,
obviously, if $C$ and $C^{\prime }$ are rationalized by preference relations 
$\succsim $ and $\trianglerighteq $, respectively, in the sense that $%
C=M(\cdot ,\succsim )$ and $C^{\prime }=M(\cdot ,\trianglerighteq ),$ then $%
d_{\text{K}}(C,C^{\prime })=D(\succsim ,\trianglerighteq ).$ On the other
hand, Eliaz and Ok \cite{E-O} have shown that every (nonempty-valued) choice
correspondence $C$ on $X$ that satisfies a slight relaxation of the
classical weak axiom of revealed preference is indeed of the form $S\mapsto
M(S,\succsim )$ for some (transitive but possibly incomplete) preference
relation $\succsim $ on $X$. It follows that we may think of $D(\succsim
,\trianglerighteq )$ as measuring the distance between $\succsim $ and $%
\trianglerighteq $ by looking at the discrepancy between the
\textquotedblleft rational choices\textquotedblright\ induced by these
preferences.

Having said this, counting the number of elements of $M(S,\succsim
)\triangle M(S,\trianglerighteq )$ is only one way of measuring the
\textquotedblleft size\textquotedblright\ of this set. Especially if there
is reason to treat the alternatives in $X$ in a non-neutral way, we may wish
to gauge this \textquotedblleft size\textquotedblright\ by means of a
measure on $2^{X}$ distinct from the counting measure.\footnote{%
Due to the political spectrum of the country, a political analyst studying
voter preferences in the case of 2020 elections may wish to weigh the
importance of the (1) Trump-Pence and (2) Biden-Harris tickets more than (3)
Hawkins-Walker and (4) Jorgensen-Cohen tickets, \textit{independently of
voter preferences}. This analyst may then choose to use a measure $\mu $
which weighs the candidates (1) and (2) more than the candidates (3) and (4)
when deciding on the size of the disagreements of the maximal sets with
respect to these preferences.} This idea yields the semimetric%
\begin{equation*}
D^{\mu }(\succsim ,\trianglerighteq )=\sum_{S\subseteq X}\mu (M(S,\succsim
)\triangle M(S,\trianglerighteq ))
\end{equation*}%
where $\mu $ is some measure on $2^{X}$. We refer to $D^{\mu }$ as the $\mu $%
\textit{-top-difference semimetric}. Obviously, $D^{\mu }=D$ where $\mu $ is
the counting measure.

We shall show later that these semimetrics act as metrics in the case of
partial orders, or complete preference relations, among other situations.%
\footnote{$D^{\mu }$ fails to distinguish between two preferences simply
because indifference and incomparability sometimes have the same effect on
maximal sets. For example, $D^{\mu }$ judges the difference between two
preferences, one exhibiting indifference everywhere and the other
incomparability everywhere, as zero. Loosely speaking, on any domain of
preferences in which indifference and incomparability are not exchangeable
(which is trivially the case if we assume away incomparabilities), each $%
D^{\mu }$ assigns a positive distance to any pair of distinct preferences.}
More important, unlike $d_{\text{KSB}},$ they are primed to evaluate the
dissimilarity of preference relations from the perspective of choice. For
instance, we have $D(\succsim ,\succsim _{1})=16>2=D(\succsim ,\succsim _{2})
$ in the case of Example 1.1, while $D(\succsim ,\succsim
_{1})=8<17=D(\succsim ,\succsim _{2})$ in the case of Example 1.2.\footnote{%
More generally, we have $D^{\mu }(\succsim ,\succsim _{1})>D^{\mu }(\succsim
,\succsim _{2})$ for every measure $\mu $ with $\mu (\{x_{1},x_{2}\})>\frac{1%
}{8}\mu (\{x_{4},x_{5}\})$ in the context of Example 1.1, while in Example
1.2, $D^{\mu }(\succsim ,\succsim _{1})>D^{\mu }(\succsim ,\succsim _{2})$
for every measure $\mu $.}

One of the main advantages of the Kemeny-Snell-Bogart metric is its
axiomatization. We thus begin our work in this paper by characterizing the
class of all $D^{\mu }$ semimetrics (where $\mu $ varies over all measures
on $2^{X}$) by means of two simple axioms. (In the case one wishes to allow
for indifferences, a third axiom is needed.) As in all axiomatizations,
these postulates allow us break down what is actually involved in measuring
the dissimilarity of two preferences by using $D^{\mu }$. We also find that
appending an additional postulate to this system, one that reflects the
neutrality of the alternatives, yields a complete characterization of the
top-difference semimetric $D,$ singling out this semimetric as a focal
element of this class.

Our axioms are built on the idea of perturbing a given preference relation
in a minimal way (so that the dissimilarity comparison is straightforward),
and then using such perturbations finitely many times to \textit{define} a
metric segment (in terms of the target semimetric). The nature of these
perturbations, and the fundamental fact that any one preference relation can
be transformed into any other given preference by applying them finitely
many times in the right order, is explained in Section 2.3, right after we
introduce the basic nomenclature of the paper. In Section 3, we formally
define our semimetrics, and show that they act as metrics in most cases of
interest. And then, in Sections 3.2 and 3.3, we introduce our axiomatic
system, and prove our characterization theorems. In Section 3.4, we show
that $D$ is the only member of the $D^{\mu }$ class which is at the same
time a weighted form of the Kemeny-Snell-Bogart metric. This highlights the
importance of $D$ even further. Finally, in Section 3.5, we obtain an
alternative formula for $D^{\mu },$ whose computation takes at most
polynomial time with respect to the size of $X$ (just like the
Kemeny-Snell-Bogart metric), and use this to obtain a very efficient method
of evaluating $D$ in the case of linear orders.

When we compute the distance between preferences by $D,$ it is difficult to
understand the significance of this magnitude (or lack thereof) without a
benchmark (while of course this quantity can always be used to make
comparisons). For this reason, in Section 4, we turn to studying the \textit{%
diameter} of certain subsets of preferences in terms of $D.$ Even for
relatively small $X$ (with about 20 elements), there are an immense number
of preference relations over $X,\footnote{%
As a side note, we note that the number of all preorders (which is the same
as that of all topologies) and the number of all partial orders (which
equals that of all $T_{0}$-topologies) on an arbitrary finite set are
presently known only up to sets with 16 elements. This is an intense area of
research in enumerative combinatorics, but the results are mainly of
asymptotic nature, as in the famous work of Kleitman and Rothschild \cite%
{Kleitmen}.}$ and this makes such diameter computations very hard.
Fortunately, however, we were able to compute this diameter exactly (Theorem
4.1) in the case of complete preferences. When $X$ is small (but still
relevant for experimental work), the resulting diameter is quite manageable
(for, say, normalization purposes).\footnote{%
If $X$ contains four elements, the largest $D$ distance between two complete
preferences is 26. For $5$-element $X$ this number goes up to $70,$ and in
the $6$-element case to $178.$ Some other computations are reported in Table
1 below.}

In Section 5, we turn to an application of our metrics $D^{\mu },$ and study
the following best approximation problem: Among all transitive extensions of
an acyclic preference relation (with no indifferences), which one is the
closest to that relation with respect to $D^{\mu }$? We find that the answer
is the transitive closure of that relation (for any $\mu $), and provide
some examples to show that this is not at all a trivial observation.

In Section 6, we turn to the problem of isometrically embedding our
preference metric spaces in a Euclidean space. When there are no
indifferences, it is known that this can be done for the Kemeny-Snell-Bogart
metrization. We prove in Section 6 that this is not possible in the case of $%
D^{\mu }$ (for any $\mu $ and any $X$ with $\left\vert X\right\vert >2$), so
one has to adopt non-Euclidean methods when working with these metrics, such
as when solving best approximation or least squares problems. The paper
concludes with a short section that points to a few avenues for future
research.

\section{Preliminaries}

\subsection{Order-Theoretic Terminology\protect\footnote{%
We summarize in this subsection all the order-theoretic concepts we use in
this paper. However, for a comprehensive treatment of these notions we
should refer the reader to authorative texts like Caspard, Leclerc and
Monjardet \cite{Caspard} and Schr\"{o}der \cite{Schroder}.}}

By a \textit{binary relation} $R$ on a nonempty set $X$, we mean any
nonempty subset of $X\times X,$ but we often adopt the usual convention of
writing $x$ $R$ $y$ instead of $(x,y)\in $ $R$. In turn, we simply write $x$ 
$R$ $y$ $R$ $z$ to mean $x$ $R$ $y$ and $y$ $R$ $z,$ and so on. The \textit{%
principal filter} and \textit{principal ideal} of any $x\in X$ with respect
to $R$ are defined as%
\begin{equation*}
x^{\downarrow ,R}:=\{a\in X:x\text{ }R\text{ }a\}\hspace{0.2in}\text{and}%
\hspace{0.2in}x^{\uparrow ,R}:=\{a\in X:a\text{ }R\text{ }x\},
\end{equation*}%
respectively. When either $x$ $R$ $y$ or $y$ $R$ $x,$ we say that $x$ and $y$
are $R$\textit{-comparable}, and put%
\begin{equation*}
\text{Inc}(R):=\{(x,y)\in X\times X:x\text{ and }y\text{ are not }R\text{%
-comparable}\}.
\end{equation*}%
If Inc$(R)=\varnothing $, we say that $R$ is \textit{total} but note that
economic theorists often refer to total relations as \textit{complete}
relations. 

For any \thinspace $S\subseteq X,$ by $x$ $R$ $S,$ we mean $x$ $R$ $y$ for
every $y\in $\thinspace $S$, and interpret the statement $S$ $R$ $x$
analogously. The set of all $R$\textit{-maximum} and $R$\textit{-maximal}
elements of $S$ are denoted by $m(S,R)$ and $M(S,R),$ respectively, that is, 
\begin{equation*}
m(S,R):=\{x\in S:x\text{ }R\text{ }S\}\text{\hspace{0.1in}and\hspace{0.1in}}%
M(S,R):=\{x\in S:y\text{ }R^{>}\text{ }x\text{ for no }y\in S\}\text{,}
\end{equation*}%
where $R^{>}$ stands for the \textit{asymmetric}\ \textit{part}\ of $R$
which is the binary relation on $X$ defined by $x$ $R^{>}$ $y$ iff $x$ $R$ $%
y $ and not $y$ $R$ $x$. (In turn, the \textit{symmetric part}\textbf{\ }of $%
R$\ is defined as $R\backslash R^{>}$.) In general, $m(S,R)\subseteq M(S,R),$
but not conversely, while $m(S,R)=M(S,R)$ whenever $R$ is total. Note also
that $M(S,R)=M(S,R^{>}).$

We denote the diagonal of $X\times X$ by $\Delta _{X},$ that is,%
\begin{equation*}
\Delta _{X}:=\{(x,x):x\in X\}\text{.}
\end{equation*}%
If $\Delta _{X}\subseteq R,$ we say that $R$ is \textit{reflexive}, and if $%
R\backslash R^{>}\subseteq \Delta _{X},$ we say that it is \textit{%
antisymmetric}. Of particular importance for the present paper is the notion
of acyclicity. We say that $R$ is \textit{acyclic} if there do not exist any
finitely many (pairwise) distinct $z_{1},...,z_{k}\in X$ such that $z_{1}$ $%
R^{>}\cdot \cdot \cdot z_{k}$ $R^{>}$ $z_{1}.$ This is a weaker property
than transitivity. Indeed, $R$ is said to be \textit{transitive} if $x$ $R$ $%
y$ $R$ $z$ implies $x$ $R$ $z,$ and \textit{quasitransitive} if $R^{>}$ is
transitive. It is plain that transitivity of a binary relation implies its
quasitransitivity, and its quasitransitivity implies its acyclicity, but not
conversely.

We say that $R$ is a \textit{preorder} on $X$ if it is reflexive and
transitive. (Total preorders are often called \textit{weak orders} in the
literature.) If, in addition, it is antisymmetric, $R$ is said to be a 
\textit{partial order} on $X,$ and if it is total, antisymmetric and
transitive, it is said to be a \textit{linear order }on $X.$ We say that $R$
is an \textit{acyclic order} (or sometimes an \textit{acyclic preference})
on $X$ if it is reflexive and acyclic. In what follows, we will denote a
generic acyclic order by $\succsim $ or $\trianglerighteq ,$ and the
asymmetric parts of these relations by $\succ $ and $\vartriangleright $,
respectively. We note that acyclic orders can always be identified with
directed acyclic graphs, which are of primary importance for many
subdisciplines of operations research.

\bigskip

\noindent \textit{Notation.} The set of all acyclic orders on $X$ is denoted
by $\mathbb{A}(X),$ that of all preorders on $X$ by $\mathbb{P}(X)$, and
that of all total preorders by $\mathbb{P}_{\text{total}}(X).$ In turn, we
denote the set of all partial orders on $X$ by $\mathbb{P}^{\ast }(X),$ and
finally, that of all linear orders on $X$ by $\mathbb{L}(X).$ Obviously,%
\begin{equation*}
\mathbb{L}(X)\subseteq \mathbb{P}^{\ast }(X)\subseteq \mathbb{P}(X)\subseteq 
\mathbb{A}(X)\hspace{0.2in}\text{and\hspace{0.2in}}\mathbb{L}(X)\subseteq 
\mathbb{P}_{\text{total}}(X)\subseteq \mathbb{P}(X).
\end{equation*}

\bigskip

Finally, we recall that the \textit{transitive closure} of a binary relation 
$R$ on $X$ is the smallest transitive relation on $X$ that contains $R$; we
denote this relation by tran$(R)$. This relation always exists; we have $x$
tran$(R)$ $y$ iff $x=x_{0}$ $R$ $x_{1}$ $R$ $\cdot \cdot \cdot $ $R$ $%
x_{k}=y $ holds for some nonnegative integer $k$ and $x_{0},...,x_{k}\in X.$
Obviously, tran$(R)$ is a preorder on $X,$ provided that $R$\ is reflexive.

\subsection{Preferences}

The standard practice of economics is to model the preference relation of an
individual as a total preorder. When one is interested in modeling the
indecisiveness of an individual over some alternatives (as in the literature
on incomplete preferences that started with Aumann \cite{aumann}), or wish
to model incomparability of some alternatives (because the outside observer
has limited data), a preference relation is taken as any preorder on the
alternative set $X.$ There are also many studies, say, in voting theory and
stable matching, where the space of preferences is identified with that of
all linear orders, or partial orders.

In all these situations, the preferences are assumed to be transitive. This
stems from focusing on \textquotedblleft rational\textquotedblright\
preferences, but on closer scrutiny, one observes that transitivity is often
a sufficient (and very convenient) property, but there are weaker
alternatives to it. For instance, one major problem with non-transitive
preferences is that these may not be maximized on some finite menus, but the
following well-known, and easily proved, fact shows that this is not a cause
for concern in the case of acyclic orders.

\bigskip

\noindent \textbf{Lemma 2.1.} \textit{Let }$X$\textit{\ be a nonempty set
and }$R$\textit{\ a reflexive binary relation on }$X.$\textit{\ Then, }$%
M(S,R)\neq \varnothing $\textit{\ for every nonempty finite }$S\subseteq X$%
\textit{\ if, and only if, }$R$\textit{\ is acyclic.}

\bigskip 

Another common rationality argument for transitivity is through the
so-called money pump arguments, but these too do not work against the
property of acyclicity. In addition, the literature on choice theory
provides plenty of rationality axioms that justify the acyclicity of \textit{%
revealed }preferences; see, among many others, \cite{Sen,J-L}. In what
follows, therefore, we model preferences on $X$ as acyclic orders on $X.$
This admits all of the standard ways of modeling preferences in economics as
special cases, and still reflect due rationality on the part of the
individuals. 

As we discussed in Section 1, our primary objective is to turn $\mathbb{A}(X)
$ into a (semi)metric space in a way that semimetric of the space reflect
the dissimilarity of two acyclic preferences on the basis of their
implications for choice. We do this in the context of a finite set of
alternatives. Thus, henceforward, we always take $X$ as a finite set that
contains at least two elements, unless otherwise is explicitly stated. (We
denote the cardinality of $X$ by $n.$) By a \textit{menu} in $X,$ we mean
any $S\subseteq X$ with $\left\vert S\right\vert \geq 2$. 

\subsection{Perturbations of Acyclic Preferences}

Let $\succsim $ be an acyclic order on $X,$ and take any distinct $a,b\in X$%
. Suppose first that $a$ and $b$ are not $\succ $-comparable. In that case
we define%
\begin{equation*}
R=\left\{ 
\begin{array}{ll}
\succsim \sqcup \{(a,b)\}, & \text{if }(a,b)\in \text{ Inc}(\succsim ), \\ 
\succsim \backslash \{(b,a)\}, & \text{if }b\sim a\text{,}%
\end{array}%
\right.
\end{equation*}%
which is a binary relation on $X$ that may or may not be acyclic.\footnote{%
For instance, where $X=\{a,b,c\},$ the binary relation $\succsim $ $%
:=\{(c,b),(b,a)\}\sqcup \triangle _{X}$ belongs to $\mathbb{A}(X),$ but $%
\succsim \oplus (a,b)$ does not.} (Here $\sim $ stands for the symmetric
part of $\succsim $.) Provided that it is acyclic, we say that $R$ is 
\textit{obtained from }$\succsim $\textit{\ by a single addition} (of $(a,b)$%
), and denote it as%
\begin{equation*}
\succsim \oplus (a,b).
\end{equation*}%
In words, $\succsim \oplus (a,b)$ is the acyclic order on $X$ that is
obtained from $\succsim $ by placing $a$ strictly above $b$ (while $\succsim 
$ itself does not render a strict ranking between $a$ and $b$). See Figures
3 and 4.

\ifx\JPicScale\undefined

\fi

\unitlength.7 mm 
\begin{picture}(65,68)(25,0)

\put(55,28){\line(0,1){20}}
\multiput(55,49)(0.12,-0.24){87}{\line(0,-1){0.12}}
\multiput(45,29)(0.12,0.24){87}{\line(0,1){0.12}}

\put(55,29){\circle*{2}}
\put(55,49){\circle*{2}}
\put(45,29){\circle*{2}}
\put(65,29){\circle*{2}}

\put(59,50){\makebox(0,0)[cc]{{\footnotesize {$x$}}}}
\put(42,29){\makebox(0,0)[cc]{{\footnotesize {$y$}}}}
\put(68,29){\makebox(0,0)[cc]{{\footnotesize {$b$}}}}
\put(59,29){\makebox(0,0)[cc]{{\footnotesize {$a$}}}}

\put(55,20){\makebox(0,0)[cc]{{\footnotesize {$\succsim$}}}}

\linethickness{0.2mm}
\multiput(100,57)(0.12,-0.24){118}{\line(0,-1){0.12}}

\multiput(93,43)(0.06,0.12){118}{\line(0,1){0.12}}

\linethickness{0.2mm}

\put(100,57){\circle*{2}}

\linethickness{0.2mm}

\put(93,43){\circle*{2}}
\put(107,43){\circle*{2}}
\put(114,29){\circle*{2}}

\put(112,43){\makebox(0,0)[cc]{{\footnotesize {$a$}}}}
\put(119,29){\makebox(0,0)[cc]{{\footnotesize {$b$}}}}

\linethickness{0.2mm}

\put(104,58){\makebox(0,0)[cc]{{\footnotesize {$x$}}}}
\put(88,43){\makebox(0,0)[cc]{{\footnotesize {$y$}}}}

\put(100,20){\makebox(0,0)[cc]{{\footnotesize {$\succsim \oplus (a,b)$}}}}

\linethickness{0.2mm}

\multiput(149,49)(0.12,-0.24){87}{\line(0,-1){0.12}}
\multiput(139,29)(0.12,0.24){87}{\line(0,1){0.12}}

\put(139,29){\circle*{2}}
\put(149,49){\circle*{2}}
\put(159,29){\circle*{2}}
\put(169,39){\circle*{2}}

\put(153,50){\makebox(0,0)[cc]{{\footnotesize {$x$}}}}
\put(135,29){\makebox(0,0)[cc]{{\footnotesize {$y$}}}}
\put(163,29){\makebox(0,0)[cc]{{\footnotesize {$a$}}}}
\put(173,39){\makebox(0,0)[cc]{{\footnotesize {$b$}}}}

\put(155,20){\makebox(0,0)[cc]{{\footnotesize {$\succsim \ominus (x,b)$}}}}

\put(104,7){\makebox(0,0)[cc]{{\bf {\footnotesize Figure 3}}}}

\end{picture}

Now suppose $a\succ b$ holds instead. Then $\succsim \backslash \{(a,b)\}$
is acyclic, and in this case we say that this relation is \textit{obtained
from }$\succsim $\textit{\ by a single deletion} (of $(a,b)$), and denote it
as 
\begin{equation*}
\succsim \ominus (a,b).
\end{equation*}%
In words, $\succsim \ominus (a,b)$ is the acyclic order on $X$ that is
obtained from $\succsim $ by eliminating the strictly higher ranking of $a$
over $b$ within $\succsim $. See Figures 3 and 4.

\ifx\JPicScale\undefined

\fi

\unitlength.7 mm 
\begin{picture}(65,68)(25,0)

\linethickness{0.2mm}
\multiput(70,57)(0.12,-0.12){118}{\line(0,-1){0.12}}
\multiput(55.96,42.86)(0.12,0.12){118}{\line(0,1){0.12}}
\multiput(70,28.72)(0.12,0.12){118}{\line(0,1){0.12}}
\multiput(55.96,42.86)(0.12,-0.12){118}{\line(0,-1){0.12}}

\put(70,57){\circle*{2}}
\put(84,43){\circle*{2}}
\put(56,43){\circle*{2}}

\put(88,42){\makebox(0,0)[cc]{{\footnotesize {$b$}}}}

\put(70,61){\makebox(0,0)[cc]{{\footnotesize {$x$}}}}
\put(52,42){\makebox(0,0)[cc]{{\footnotesize {$a$}}}}

\put(70,18){\makebox(0,0)[cc]{{\footnotesize {$\succsim$}}}}

\put(70,29){\circle*{2}}
\put(70,25){\makebox(0,0)[cc]{{\footnotesize {$y$}}}}

\linethickness{0.2mm}
\put(106,29){\line(0,1){28}}
\put(106,57){\circle*{2}}
\put(106,29){\circle*{2}}
\put(106,48){\circle*{2}}
\put(106,38.5){\circle*{2}}

\put(110,58){\makebox(0,0)[cc]{{\footnotesize {$x$}}}}

\put(110,48){\makebox(0,0)[cc]{{\footnotesize $a$}}}
\put(110,38.5){\makebox(0,0)[cc]{{\footnotesize $b$}}}
\put(110,28){\makebox(0,0)[cc]{{\footnotesize $y$}}}

\linethickness{0.2mm}

\multiput(130,49)(0.12,-0.12){168}{\line(0,1){0.12}}
\multiput(150,28.72)(0.12,0.12){118}{\line(0,1){0.12}}

\put(140,39){\circle*{2}}
\put(130,49){\circle*{2}}
\put(164,43){\circle*{2}}
\put(150,29){\circle*{2}}

\put(126,49){\makebox(0,0)[cc]{{\footnotesize {$x$}}}}
\put(137,36){\makebox(0,0)[cc]{{\footnotesize {$a$}}}}
\put(168,43){\makebox(0,0)[cc]{{\footnotesize {$b$}}}}
\put(150,25){\makebox(0,0)[cc]{{\footnotesize {$y$}}}}

\put(108,18){\makebox(0,0)[cc]{{\footnotesize {$\succsim \oplus (a,b)$}}}}

\put(151,18){\makebox(0,0)[cc]{{\footnotesize {$\succsim \ominus (x,b)$}}}}

\put(108,7){\makebox(0,0)[cc]{{\bf {\footnotesize Figure 4}}}}

\end{picture}

We emphasize that both $\succsim \oplus (a,b)$ and $\succsim \ominus (a,b)$
belong to $\mathbb{A}(X)$. (In the first case this is true by definition,
and in the second case this is true by necessity.) Moreover, when $a$ and $b$
are not $\succsim $-comparable, we have $(\succsim \oplus (a,b))\ominus
(a,b)=$ $\succsim ,$ and similarly, when $a\succ b$, we have $(\succsim
\ominus (a,b))\oplus (a,b)=$ $\succsim $. However, when $a\sim b,$ we have $%
(\succsim \oplus (a,b))\ominus (a,b)=$ $\succsim \backslash \{(a,b),(b,a)\}$.

Let $\succsim _{0}$ and $\trianglerighteq $ be two acyclic orders on $X.$ We
say that $\succsim _{0}$ is a \textit{one-step perturbation} \textit{of} $%
\succsim $ \textit{toward }$\trianglerighteq $ if either (i) 
\begin{equation}
\succsim _{0}\text{\thinspace }=\text{ }\succsim \ominus (a,b)\hspace{0.2in}%
\text{and\hspace{0.2in}not }a\vartriangleright b  \label{per2}
\end{equation}%
and%
\begin{equation}
x\succ b\hspace{0.2in}\text{for every }x\in X\text{ with }x\vartriangleright
b  \label{proper}
\end{equation}%
for some $(a,b)\in $ $\succ $; or (ii)%
\begin{equation}
\succsim _{0}\text{\thinspace }=\text{ }\succsim \oplus (a,b)\hspace{0.2in}%
\text{and\hspace{0.2in}}a\vartriangleright b  \label{per1}
\end{equation}%
for some $(a,b)\in $ Inc$(\succ )$. Intuitively speaking, when this is the
case, we understand that the ranking positions of $a$ and $b$ in $\succsim $
is altered in a way that becomes identical to how these elements are ranked
by $\trianglerighteq $. (This is captured by (\ref{per2}) and (\ref{per1}).)
In this sense, we think of $\succsim _{0}$ as \textquotedblleft more
similar\textquotedblright\ to $\trianglerighteq $ than $\succsim $ is. This
viewpoint is further enforced by the requirement (\ref{proper}) which
maintains that the ordering of $b$ in $\succsim $ is consistent with that in 
$\trianglerighteq $. The following example highlights the importance of this
consistency condition.

\bigskip

\noindent \textit{Example 2.1.} Let $X=\{a,b,c\},$ and consider the acyclic
orders $\succsim $ and $\trianglerighteq $ on $X$ with $\succsim $ $:=$%
\thinspace $\Delta _{X}\sqcup \{(a,b)\}$ and $\trianglerighteq $ $:=$%
\thinspace $\Delta _{X}\sqcup \{(c,b)\}.$ Then,\thinspace $\Delta _{X}=$ $%
\succsim \ominus (a,b),$ but $\Delta _{X}$ is not a one-step perturbation of 
$\succsim $ toward $\trianglerighteq $. Indeed, in this case, it is not
really evident whether or not $\Delta _{X}$ is \textquotedblleft more
similar\textquotedblright\ to $\trianglerighteq $ than $\succsim $ is,
especially if we focus on the maximal elements in various subsets of $X.$ If
we restrict attention to the sets $\{a,c\}$ and $\{c,b\},$ the behavior of $%
\succsim $ and $\Delta _{X}$ are identical, while on $\{a,b\}$ the behavior
of $\Delta _{X}$ is identical to that of $\trianglerighteq $. However, on
the grand set $X,$ the diagonal relation $\Delta _{X}$ behaves quite
differently than $\trianglerighteq $. Indeed, \thinspace $\Delta _{X}$
declares $b$ as maximal in $X,$ while $b$ is minimal in $X$ relative to $%
\trianglerighteq $. By contrast, $\succsim $ and $\trianglerighteq $ have
the same set of maximal elements in $X.$ We impose the consistency condition
(\ref{proper}) on one-step perturbations precisely to avoid such ambiguous
situations. $\Vert $

\bigskip

In what follows, if $\succsim _{0}$ is a one-step perturbation of $\succsim $
toward $\trianglerighteq ,$ we write%
\begin{equation*}
\succsim \text{\thinspace }\rightarrow \text{\thinspace }\succsim _{0}\text{%
\thinspace }\twoheadrightarrow \text{\thinspace }\trianglerighteq \text{.}
\end{equation*}%
Generalizing this concept, for any integer $n\geq 2,$ we say an acyclic
order $\succsim _{n-1}$ on $X$ is an $n$\textit{-step perturbation of }$%
\succsim $\textit{\ toward} $\trianglerighteq ,$ if there exist $\succsim
_{0},...,\succsim _{n-2}$ $\in \mathbb{A}(X)$ such that $\succsim $%
\thinspace $\rightarrow $\thinspace $\succsim _{0}$\thinspace $%
\twoheadrightarrow $\thinspace $\trianglerighteq $ and%
\begin{equation*}
\succsim _{k-1}\text{\thinspace }\rightarrow \text{\thinspace }\succsim _{k}%
\text{\thinspace }\twoheadrightarrow \text{\thinspace }\trianglerighteq 
\text{ for each }k=1,...,n-1.
\end{equation*}%
Finally, we say that an $\succsim _{\ast }\in \mathbb{A}(X)$ is \textit{%
in-between} $\succsim $ and $\trianglerighteq $ if $\succsim _{\ast }$ is an 
$n$\textit{-step perturbation of }$\succsim $\textit{\ toward} $%
\trianglerighteq $ for some positive integer $n.$ And if $\succsim _{\ast }$ 
$=$ $\trianglerighteq $ here, we say that $\succsim $ is \textit{transformed
into} $\trianglerighteq $ \textit{in finitely many steps.}

\bigskip

\noindent \textit{Remark.} In the literature on metrics on preference
relations, one often says that a binary relation $R_{0}$ on $X$ is
\textquotedblleft between\textquotedblright\ the binary relations $R_{\ast }$
and $R^{\ast }$ on $X$ if $R_{\ast }\cap R^{\ast }\subseteq R\subseteq
R_{\ast }\cup R^{\ast }$.(See, for instance, \cite{Bogart1,Bogart2,Cook2}.)
Our definition of being \textquotedblleft in-between\textquotedblright\ is
more stringent than this concept, due to the consistency condition (\ref%
{proper}). For instance, in the context of Example 2.1, $\Delta _{X}$ is
\textquotedblleft between\textquotedblright\ $\succsim $ and $%
\trianglerighteq $ according to the betweenness definition of the
literature, but $\Delta _{X}$ is not in-between\ $\succsim $ and $%
\trianglerighteq $ according to our definition. This is consistent with the
main motivation of the present work. We would like to think of an acyclic
order $\succsim _{\ast }$ on $X$ that is in-between $\succsim $ and $%
\trianglerighteq $ as one that is \textquotedblleft more
similar\textquotedblright\ in its order structure to $\trianglerighteq $
than $\succsim $ is. As we have seen in Example 2.1, at least insofar as
which elements are declared maximal in various menus, being
\textquotedblleft between\textquotedblright\ two acyclic orders does not
fully support this interpretation.

\bigskip

The following result provides the fundamental force behind the
axiomatization that we present in the next section.

\bigskip

\noindent \textbf{Theorem 2.2.} \textit{Let }$\succsim $ \textit{and }$%
\trianglerighteq $ \textit{be distinct} \textit{acyclic orders on }$X$ 
\textit{with the same symmetric parts.} \textit{Then,} $\succsim $ \textit{%
can be} \textit{transformed into} $\trianglerighteq $ \textit{by finitely
many one-step perturbations}.\footnote{%
Example 2.1 points to the nontriviality of this claim. Arbitrary addition
and/or deletions of pairs of alternatives from $\succsim $ may not be able
to transform $\succsim $ into $\trianglerighteq $. Instead, the theorem
claims that there is always a \textquotedblleft right\textquotedblright\
order of doing these perturbations which would transform $\succsim $ into $%
\trianglerighteq $.}

\begin{proof}
We will prove that there exists an $\succsim _{0}$ $\in \mathbb{A}(X)$ such
that $\succsim $\thinspace $\rightarrow $\thinspace $\succsim _{0}$%
\thinspace $\twoheadrightarrow $\thinspace $\trianglerighteq $. The more
general statement of the theorem will then follow by induction.

Note first that if $\succ $ $\subseteq $ $\vartriangleright ,$ then the
containment is proper (because $\succsim $ $\neq $\textit{\ }$%
\trianglerighteq $), so we are readily done by setting $\succsim _{0}$%
\thinspace $:=$ $\succsim \oplus (a,b)$ for any $a,b\in $ $\vartriangleright
\backslash \succ $. Let us then assume that $\succ $ is not contained within 
$\vartriangleright $, that is, 
\begin{equation*}
B:=\{b\in X:(a,b)\in \,\succ \backslash \vartriangleright \text{ for some }%
a\in X\}\neq \varnothing \text{.}
\end{equation*}%
We pick any tran$(\succ )$-minimal element $b^{\ast }$ of $B,$ and any $%
a^{\ast }\in X$ with $a^{\ast }\succ b^{\ast }$ but not $a^{\ast
}\vartriangleright b^{\ast }$. If 
\begin{equation*}
x\succ b^{\ast }\hspace{0.2in}\text{for every }x\in X\text{ with }%
x\vartriangleright b^{\ast },
\end{equation*}%
then we are done by setting $\succsim _{0}$ $:=$ $\succsim \ominus (a,b).$
We thus assume that this is not the case, that is, there is an $x\in X$ such
that%
\begin{equation}
x\vartriangleright b^{\ast }\hspace{0.2in}\text{and}\hspace{0.2in}\text{not }%
x\succ b^{\ast }\text{.}  \label{x}
\end{equation}

Next, we define $\succsim _{0}$ $:=$ $\succsim \sqcup \{(x,b^{\ast })\}.$
Given that $x\vartriangleright b^{\ast },$ our proof will be complete if we
can show that $\succsim _{0}$ $=$ $\succsim \oplus (x,b^{\ast }).$ But note
that we cannot have $b^{\ast }\succ x$ here, because otherwise $x\in B,$ and 
$b^{\ast }\succ x$ contradicts the tran$(\succ )$-minimality of $b^{\ast }$
in $B$. We cannot have $b^{\ast }\sim x$ either, because $x\vartriangleright
b^{\ast }$ while $\sim $ equals to the symmetric part of $\trianglerighteq $
by hypothesis. Thus: $(x,b^{\ast })\in $ Inc$(\succsim ).$ By definition of
the relation $\succsim \oplus (x,b^{\ast }),$ it thus remains only to show
that $\succsim _{0}$ is acyclic. To derive a contradiction, suppose this is
not the case, that is, assume there exist an $n\in \mathbb{N}$ and distinct $%
z_{1},...,z_{n}\in X$ with $z_{1}\succ _{0}\cdot \cdot \cdot \succ
_{0}z_{n}\succ _{0}z_{1}.$ Since $\succsim $ $\in $ $\mathbb{A}(X),$ we must
have $(z_{k},z_{k+1(\func{mod}n)})=(x,b^{\ast })$ for some $k=1,...,n.$
Thus, relabelling if necessary, we may assume that $(z_{1},z_{2})=(x,b^{\ast
})$ in which case we have 
\begin{equation}
b^{\ast }=z_{2}\succ \cdot \cdot \cdot \succ z_{n}\succ z_{1}  \label{bb}
\end{equation}%
by definition of $\succsim _{0}$. Now, if $z_{k}\vartriangleright z_{k+1(%
\func{mod}n)}$ for each $k=2,...,n,$ then%
\begin{equation*}
b^{\ast }=z_{2}\vartriangleright \cdot \cdot \cdot \vartriangleright
z_{n}\vartriangleright z_{1}=x\vartriangleright b^{\ast }
\end{equation*}%
and we contradict the acyclicity of $\trianglerighteq $. Let us then assume
that $z_{k}\vartriangleright z_{k+1(\func{mod}n)}$ fails for some $k=2,...,n$%
. In view of (\ref{bb}), this means that $z_{k}\in B$ for some $k\in
\{1,...,n\}\backslash \{2\}.$ But again by (\ref{bb}), we have $b^{\ast }$
tran$(\succ )$ $z_{k}$ for every $k\in \{1,...,n\}\backslash \{2\},$ so this
finding contradicts the tran$(\succ )$-minimality of $b^{\ast }$ in $B.$ We
conclude that $\succsim _{0}$ $\in $ $\mathbb{A}(X).$ As noted above, this
completes the proof.
\end{proof}

In Figure 5, we provide a simple illustration of how a partial order (in
this case the pentagon lattice) is transformed into another by means of
three one-step perturbations. In this example, the middle two partial orders
are in-between left-most and right-most partial orders. (In particular, we
have $\succsim ^{\ast }$ $=$ $((\succsim \oplus (y,z))\ominus (w,a))\ominus
(z,a)$.) But despite what this example may suggest, we emphasize that a
non-transitive (but always acyclic) binary relation may be in-between two
partial orders.

\ifx.7\undefined

\fi

\unitlength0.8 mm 
\begin{picture}(65,68)(25,0)

\linethickness{0.2mm}
\put(48,30){\line(0,1){10}}
\multiput(38,50)(0.12,-0.12){83}{\line(1,0){0.12}}
\multiput(38,20)(0.12,0.12){83}{\line(1,0){0.12}}
\multiput(28,35)(0.12,-0.18){83}{\line(0,-1){0.18}}
\multiput(28,35)(0.12,0.18){83}{\line(0,1){0.18}}

\put(38,50){\circle*{2}}
\put(48,40){\circle*{2}}
\put(48,40){\circle*{2}}
\put(48,30){\circle*{2}}
\put(38,20){\circle*{2}}
\put(28,35){\circle*{2}}

\put(52,30){\makebox(0,0)[cc]{{\footnotesize {$w$}}}}
\put(52,40){\makebox(0,0)[cc]{{\footnotesize {$y$}}}}
\put(24,35){\makebox(0,0)[cc]{{\footnotesize {$z$}}}}
\put(38,54){\makebox(0,0)[cc]{{\footnotesize {$x$}}}}
\put(42,20){\makebox(0,0)[cc]{{\footnotesize {$a$}}}}

\put(62,50){\makebox(0,0)[cc]{{\footnotesize {$\longrightarrow$}}}}
\put(104,50){\makebox(0,0)[cc]{{\footnotesize {$\longrightarrow$}}}}
\put(149,50){\makebox(0,0)[cc]{{\footnotesize {$\longrightarrow$}}}}

\put(38,12){\makebox(0,0)[cc]{{\footnotesize {$\succsim$}}}}
\put(166,12){\makebox(0,0)[cc]{{\footnotesize {$\succsim^*$}}}}

\put(82,45){\line(0,1){10}}
\multiput(82,45)(0.12,-0.12){83}{\line(1,0){0.12}}
\multiput(72,35)(0.12,0.12){83}{\line(1,0){0.12}}
\multiput(82,25)(0.12,0.12){83}{\line(1,0){0.12}}
\multiput(72,35)(0.12,-0.12){83}{\line(1,0){0.12}}

\put(82,45){\circle*{2}}
\put(72,35){\circle*{2}}
\put(92,35){\circle*{2}}
\put(82,55){\circle*{2}}
\put(82,25){\circle*{2}}

\put(96,35){\makebox(0,0)[cc]{{\footnotesize {$w$}}}}
\put(87,45){\makebox(0,0)[cc]{{\footnotesize {$y$}}}}
\put(87,55){\makebox(0,0)[cc]{{\footnotesize {$x$}}}}
\put(68,35){\makebox(0,0)[cc]{{\footnotesize {$z$}}}}
\put(82,21){\makebox(0,0)[cc]{{\footnotesize {$a$}}}}

\put(126,45){\line(0,1){10}}
\multiput(126,45)(0.12,-0.12){83}{\line(1,0){0.12}}
\multiput(116,35)(0.12,0.12){83}{\line(1,0){0.12}}
\multiput(106,25)(0.12,0.12){83}{\line(1,0){0.12}}

\put(126,55){\circle*{2}}
\put(126,45){\circle*{2}}
\put(106,25){\circle*{2}}
\put(116,35){\circle*{2}}
\put(136,35){\circle*{2}}

\put(130,55){\makebox(0,0)[cc]{{\footnotesize {$x$}}}}
\put(130,45){\makebox(0,0)[cc]{{\footnotesize {$y$}}}}
\put(140,35){\makebox(0,0)[cc]{{\footnotesize {$w$}}}}
\put(120,35){\makebox(0,0)[cc]{{\footnotesize {$z$}}}}
\put(110,25){\makebox(0,0)[cc]{{\footnotesize {$a$}}}}

\put(166,26){\line(0,1){24}}
\multiput(166,39)(0.12,-0.12){103}{\line(1,0){0.12}}

\multiput(153,26)(0.12,0.12){103}{\line(1,0){0.12}}

\put(166,39){\circle*{2}}
\put(166,26){\circle*{2}}
\put(153,26){\circle*{2}}
\put(179,26){\circle*{2}}
\put(166,51){\circle*{2}}

\put(166,55){\makebox(0,0)[cc]{{\footnotesize {$x$}}}}
\put(169,40){\makebox(0,0)[cc]{{\footnotesize {$y$}}}}
\put(153,23){\makebox(0,0)[cc]{{\footnotesize {$w$}}}}
\put(179,23){\makebox(0,0)[cc]{{\footnotesize {$z$}}}}
\put(166,23){\makebox(0,0)[cc]{{\footnotesize {$a$}}}}

\put(100,6){\makebox(0,0)[cc]{{\bf {\footnotesize Figure 5}}}}

\end{picture}

\section{A Class of Dissimilarity Semimetrics for Preferences}

\subsection{Top-Difference Semimetrics}

For any positive measure $\mu $ on $2^{X},\footnote{%
We consider the zero measure $S\mapsto 0$ on $2^{X}$ as a member of the
family of all positive measures on $2^{X},$ but, of course, the semimetric $%
D^{\mu }$ where $\mu $ is the zero measure is of no interest.}$ we define
the $\mu $\textbf{-top-difference semimetric} $D^{\mu }:\mathbb{A}(X)\times 
\mathbb{A}(X)\rightarrow \lbrack 0,\infty )$ by 
\begin{equation*}
D^{\mu }(\succsim ,\trianglerighteq ):=\sum_{S\subseteq X}\mu (M(S,\succsim
)\triangle M(S,\trianglerighteq ))\text{.}
\end{equation*}%
In the special case where $\mu $ is the counting metric, we refer to $D^{\mu
}$ simply as the \textbf{top-difference semimetric}, and denote it by $D,$
that is, 
\begin{equation}
D(\succsim ,\trianglerighteq ):=\sum_{S\subseteq X}\left\vert M(S,\succsim
)\triangle M(S,\trianglerighteq )\right\vert   \label{old}
\end{equation}%
for any $\succsim ,\trianglerighteq $ $\in \mathbb{A}(X).$

That each $D^{\mu }$ is indeed a semimetric on $\mathbb{A}(X)$ is
straightforward. Unless $X$ is a singleton, however, $D^{\mu }$ does not act
as a metric even on $\mathbb{P}(X).$ For instance, $D^{\mu }$ cannot
distinguish between complete indifference and complete incomparability, that
is, $D^{\mu }(\Delta _{X},X\times X)=0$ for any measure $\mu $ on $2^{X}$
while $\Delta _{X}$ and $X\times X$ are distinct preorders on $X$ when $%
\left\vert X\right\vert \geq 2.$ (This is simply because the maximal
elements relative to these relations are the same in every menu.) For
another example, note that a partial order and a preorder on $X$ may have
the same asymmetric part, but may nevertheless be distinct relations on $X.$

In passing, we note that there are interesting subclasses of acyclic orders
on which $D^{\mu }$ acts as a metric, provided that $\mu $ has full support.
We present two examples to illustrate.

\bigskip

\noindent \textit{Example 3.1.} Any $D^{\mu }$ acts as a metric on the set
of all partial orders on $X.$ That is, $D^{\mu }|_{\mathbb{P}^{\ast
}(X)\times \mathbb{P}^{\ast }(X)}$ is a metric on $\mathbb{P}^{\ast }(X)$
for any measure $\mu $ on $2^{X}$. $\Vert $

\bigskip

\noindent \textit{Example 3.2.} For any preorder $\succsim $ on $X,$ the 
\textit{indifference part} of $\succsim $, denoted by ind$(\succsim ),$ is
the binary relation on $X$ defined by $(x,y)\in $ ind$(\succsim )$ iff%
\begin{equation*}
x\succ z\text{ iff }y\succ z\hspace{0.2in}\text{and\hspace{0.2in}}z\succ x%
\text{ iff }z\succ y
\end{equation*}%
for every $z\in X.$ (If we interpret $\succsim $ as the preference relation
of a person, then $(x,y)\in $ ind$(\succsim )$ means that this individual
treats $x$ and $y$ as identical objects in every menu; see \cite{E-O} and 
\cite{RR}.)

It is immediate from this definition that ind$(\succsim )$ is an equivalence
relation on $X$, and $\sim $ $\subseteq $ ind$(\succsim ).$ If $\succsim $
is total, then this holds as an equality, but in general, it may well hold
properly.\footnote{%
For instance, let $X$ consist of the 2-vectors $x=(0,5),$ $y=(5,0)$ and $%
z=(6,1),$ and let $\succsim $ be the coordinatewise ordering on $X.$ Then,
ind$(\succsim )$ contains all elements of $X\times X$ except $(y,z)$ and $%
(z,y),$ while $\sim $ equals $\Delta _{X}.$} Those preorders whose symmetric
parts match their indifference parts exactly are of immediate interest for
decision theory. Eliaz and Ok \cite{E-O} refer to a preorder $\succsim $ on $%
X$ with this property, that is, when $\sim $ $=$ ind$(\succsim ),$ as a 
\textit{regular} preorder on $X$.

Let $\succsim _{1}$ and $\succsim _{2}$\ be two regular preorders on $X$\
such that $M(S,\succsim _{1})=M(S,\succsim _{2})$ for every doubleton $%
S\subseteq X.$ We claim that $\succsim _{1}$ $=$ $\succsim _{2}$. To see
this, note that for any distinct $x,y\in X,$ we have $x\succ _{i}y$ iff $%
\{x\}=$ $M(\{x,y\},\succsim _{i})$ for $i=1,2.$ By hypothesis, therefore, $%
\succ _{1}=$ $\succ _{2}$. But then, by definition of ind$(\cdot )$, we have
ind$(\succsim _{1})=$ ind$(\succsim _{2})$ as well. Since both $\succsim
_{1} $ and $\succsim _{2}$ are regular, it follows that $\sim _{1}=$ $\sim
_{2}$.

As an immediate consequence of this observation, we see that the restriction
of $D^{\mu }$ to the class of all regular preorders on $X$ yields a metric
on that class, for any measure $\mu $ on $2^{X}$. In particular, each $%
D^{\mu }$ is a metric on the set $\mathbb{P}_{\text{total}}(X)$ of all
complete preorders, the standard setup of economic theory. $\Vert $

\subsection{Axioms}

Let $d$ be a semimetric on $\mathbb{A}(X).$ The first axiom we impose on $d$
says simply that if an acyclic order is in-between two acyclic orders on $X,$
say, $\succsim $ and\textit{\ }$\trianglerighteq ,$ then that order must lie
on the metric segment between $\succsim $ and\textit{\ }$\trianglerighteq $
relative to $d.$ That is:

\bigskip

\noindent \textbf{Axiom 1.} \textit{For any }$\succsim $, $\succsim _{0}$ 
\textit{and }$\trianglerighteq $\textit{\ in} $\mathbb{A(}X)$ \textit{such
that }$\succsim _{0}$\textit{\ is in-between }$\succsim $\textit{\ and }$%
\trianglerighteq ,$\textit{\ we have}%
\begin{equation*}
d(\succsim ,\trianglerighteq )=d(\succsim ,\succsim _{0})+d(\succsim
_{0},\trianglerighteq ).
\end{equation*}

\bigskip 

We may, of course, equivalently state this axiom in the following way which
is easier to check:

\bigskip

\noindent \textbf{Axiom 1'.} \textit{For any }$\succsim $, $\succsim _{0}$, 
\textit{and }$\trianglerighteq $\textit{\ in} $\mathbb{A(}X)$ \textit{such
that }$\succsim $\thinspace $\rightarrow $\thinspace $\succsim _{0}$%
\thinspace $\twoheadrightarrow $\thinspace $\trianglerighteq ,$\textit{\ we
have}%
\begin{equation*}
d(\succsim ,\trianglerighteq )=d(\succsim ,\succsim _{0})+d(\succsim
_{0},\trianglerighteq ).
\end{equation*}

\bigskip 

One may view these (equivalent) axioms as \textit{additivity }properties.
For instance, when $\succsim $\thinspace $\rightarrow $\thinspace $\succsim
_{0}$\thinspace $\twoheadrightarrow $\thinspace $\trianglerighteq ,$ we know
that $\succsim _{0}$ and $\trianglerighteq $ are \textquotedblleft more
similar\textquotedblright\ than $\succsim $ and $\trianglerighteq $ are, so
a metric $d$ that captures the dissimilarity of acyclic orders should
certainly declare that $d(\succsim ,\trianglerighteq )>d(\succsim
_{0},\trianglerighteq ).$ Axiom 1' says further that the \textquotedblleft
excess dissimilarity\textquotedblright\ of $\succsim $ and $\trianglerighteq 
$ additively decomposes into the dissimilarity of $\succsim $ and $\succsim
_{0}$ and that of $\succsim _{0}$ and $\trianglerighteq .$ As such, Axiom 1'
(hence Axiom 1) are not only duly compatible with how we view the notion of
one-step perturbations (and hence the concept of being in-between), but it
also brings a mathematically convenient structure for accounting the effects
of such perturbations.\footnote{%
There are many papers in the literature on metrics for preference relations
in which such additivity axioms are used; see, for instance, \cite%
{Bogart1,Bogart2,Cook2}. The difference of Axiom 1 from its predecessors
lies in the way we defined the notion of one-step perturbations, and hence
the concept of being in-between.}

We next consider two particularly simple partial orders on $X.$ For any
distinct $a,b\in X,$ we define $\succsim _{ab}$ and $\succsim _{ab}^{+}$ as
the partial orders on $X$ whose asymmetric parts are given as%
\begin{equation*}
\succ _{ab}\text{ }:=(X\backslash \{a,b\})\times \{a,b\}
\end{equation*}%
and%
\begin{equation*}
\succ _{ab}^{+}\text{ }:=\text{ }\succ _{ab}\sqcup \{(a,b)\}\text{.}
\end{equation*}%
In words, $\succsim _{ab}$ ranks every alternative other than $a$ and $b$
strictly above both $a$ and $b,$ making no other pairwise comparisons
(including that between $a$ and $b$). In turn, $\succsim _{ab}^{+}$ is the
same relation as $\succsim _{ab}$ except that it ranks $a$ strictly higher
than $b.$ (See Figure 6 for the Hasse diagrams of these partial orders in
the case where $X$ has six elements.)

\ifx\JPicScale\undefined

\fi

\unitlength.7 mm 
\begin{picture}(65,68)(25,0)

\multiput(81,51)(0.12,-0.66){36}{\line(0,-1){0.12}}
\multiput(69,52)(0.12,-0.17){140}{\line(0,-1){0.12}}

\multiput(65,29)(0.12,0.17){140}{\line(0,-1){0.12}}

\multiput(65,29)(0.12,0.66){36}{\line(0,1){0.12}}

\multiput(57,52)(0.12,-0.36){65}{\line(0,1){0.12}}
\multiput(57,52)(0.12,-0.10){235}{\line(0,1){0.12}}

\multiput(85,28)(0.12,0.36){65}{\line(0,1){0.12}}
\multiput(65,29)(0.12,0.10){235}{\line(0,1){0.12}}

\put(85,29){\circle*{2}}
\put(65,29){\circle*{2}}

\put(65,25){\makebox(0,0)[cc]{{\footnotesize {$a$}}}}
\put(85,25){\makebox(0,0)[cc]{{\footnotesize {$b$}}}}

\put(75,18){\makebox(0,0)[cc]{{\footnotesize {$\succsim_{ab}$}}}}

\put(139,38){\makebox(0,0)[cc]{{\footnotesize {$a$}}}}
\put(139,28){\makebox(0,0)[cc]{{\footnotesize {$b$}}}}

\put(137,18){\makebox(0,0)[cc]{{\footnotesize {$\succsim_{ab}^+$}}}}

\put(57,52){\circle*{2}}
\put(69,52){\circle*{2}}
\put(81,52){\circle*{2}}
\put(93,52){\circle*{2}}

\put(57,56){\makebox(0,0)[cc]{{\footnotesize {$x_1$}}}}
\put(69,56){\makebox(0,0)[cc]{{\footnotesize {$x_2$}}}}
\put(81,56){\makebox(0,0)[cc]{{\footnotesize {$x_3$}}}}
\put(94,56){\makebox(0,0)[cc]{{\footnotesize {$x_4$}}}}

\put(117,60){\circle*{2}}
\put(129,60){\circle*{2}}
\put(141,60){\circle*{2}}
\put(153,60){\circle*{2}}
\put(135,39){\circle*{2}}
\put(135,28){\circle*{2}}

\put(117,64){\makebox(0,0)[cc]{{\footnotesize {$x_1$}}}}
\put(129,64){\makebox(0,0)[cc]{{\footnotesize {$x_2$}}}}
\put(141,64){\makebox(0,0)[cc]{{\footnotesize {$x_3$}}}}
\put(154,64){\makebox(0,0)[cc]{{\footnotesize {$x_4$}}}}

\put(135,29){\line(0,1){10}}
\multiput(117,60)(0.12,-0.14){146}{\line(0,-1){0.12}}
\multiput(135,39)(0.12,0.14){146}{\line(0,-1){0.12}}

\multiput(129,60)(0.12,-0.42){52}{\line(0,-1){0.12}}
\multiput(135,39)(0.12,0.42){52}{\line(0,-1){0.12}}

\put(104,7){\makebox(0,0)[cc]{{\bf {\footnotesize Figure 6}}}}

\end{picture}

The following axiom is a neutrality property that posits that the distance
between $\succsim _{ab}$ and $\succsim _{ab}^{+}$ is independent of both $a$
and $b$, and normalizes this distance to 1.

\bigskip

\noindent \textbf{Axiom 2.} $d(\succsim _{ab},\succsim _{ab}^{+})=1$ \textit{%
for every distinct }$a,b\in X$.

\bigskip

To state our next axiom, we define%
\begin{equation*}
N(b,\succsim ):=\left\vert \left\{ x\in X\backslash \{b\}:\text{ not }x\succ
b\right\} \right\vert
\end{equation*}%
for any $b\in X$ and $\succsim $ $\in \mathbb{A}(X).$ Thus $N(b,\succsim )$
is the number of elements of $X\backslash \{b\}$ that are not ranked
strictly higher than $b$ by $\succsim $.

To understand the significance of this number, take any $\succsim $ $\in 
\mathbb{A}(X)$ and any $a,b\in X$ with $a\succ b.$ Put $\succsim _{0}$ $:=$ $%
\succsim $ $\ominus (a,b)$. Then, there are menus $S$ for which $b$ is $%
\succsim _{0}$-maximal -- that is, $b$ is a \textquotedblleft
choice\textquotedblright\ from $S$ for an individual with preferences $%
\succsim _{0}$ -- but it is not $\succsim $-maximal. This happens precisely
for those $S\subseteq X$ such that%
\begin{equation}
S=\{a,b\}\sqcup T\hspace{0.2in}\text{for some }T\subseteq N(b,\succsim ).
\label{num}
\end{equation}%
Moreover, on each such menu, the set of \textquotedblleft
choices\textquotedblright\ on the basis of $\succsim $ and $\succsim _{0}$
differ from each other by $\{b\}$ just as the set of \textquotedblleft
choices\textquotedblright\ on the basis of $\succsim _{ab}$ and $\succsim
_{ab}^{+}$ differ from each other by $\{b\}.$ Consequently, per such menu,
it makes sense to deem the dissimilarity between $\succsim $ and $\succsim
_{0}$ as the same as that between $\succsim _{ab}$ and $\succsim _{ab}^{+},$
at least insofar as we wish to capture the dissimilarity of preference
relations on the basis of what they declare maximal in various menus. As
there are $2^{N(b,\succsim )}$ many menus that satisfy (\ref{num}) (by
definition of $N(b,\succsim )),$ therefore, a consistent assignment of a
\textquotedblleft distance\textquotedblright\ between $\succsim $ and $%
\succsim _{0}$ would be $2^{N(b,\succsim )}d(\succsim _{ab},\succsim
_{ab}^{+}).$

We can reason analogously when $(a,b)\in $ Inc$(\succ )$ and $\succsim _{0}$
equals, instead, $\succsim $ $\oplus (a,b).$ In this case, a pivotal menu $S$
would be a subset of $X$ such that $a\in S$ and $b\in M(S,\succsim ).$ This
happens for those $S\subseteq X$ such that%
\begin{equation*}
S=\{b\}\sqcup T\hspace{0.2in}\text{for some }T\subseteq N(b,\succsim )\text{
with }a\in T.
\end{equation*}%
By definition of $N(b,\succsim )$ there are exactly $2^{N(b,\succsim )-1}$
many such menus, so reasoning as in the previous paragraph, we arrive at the
conclusion that a consistent assignment of a \textquotedblleft
distance\textquotedblright\ between $\succsim $ and $\succsim _{0}$ is $%
2^{N(b,\succsim )-1}d(\succsim _{ab},\succsim _{ab}^{+}).$

These considerations prompt:

\bigskip

\noindent \textbf{Axiom 3.} \textit{For any }$\succsim $ $\in \mathbb{A}(X)$%
\textit{\ and }$a,b\in X$\textit{, if }$a$\textit{\ and }$b$\textit{\ are
not }$\succ $\textit{-comparable,} 
\begin{equation*}
d(\succsim ,\succsim \oplus (a,b))=2^{N(b,\succsim )-1}d(\succsim
_{ab},\succsim _{ab}^{+}),
\end{equation*}%
\textit{and if} $a\succ b,$%
\begin{equation*}
d(\succsim ,\succsim \ominus (a,b))=2^{N(b,\succsim )}d(\succsim
_{ab},\succsim _{ab}^{+})\text{.}
\end{equation*}

\bigskip

Our final axiom is very basic. The notion of \textquotedblleft
dissimilarity\textquotedblright\ for preferences (acyclic orders) that we
focus on in this paper stems from the dissimilarity of the sets of choices
that these preferences induce on menus (subsets of $X$). And, as usual, we
model all potential choices of an individual with a given preference
relation on a menu $S$ as the set of all maximal elements of $S$ relative to
that preference. But maximal elements of a set with respect to a binary
relation depends only on the asymmetric part of that relation. That is, the
maximal subsets of any $S\subseteq X$ relative to two acyclic orders on $X$
with the same asymmetric part are identical. Thus:

\bigskip

\noindent \textbf{Axiom 4.} \textit{For any }$\succsim ,\trianglerighteq $ $%
\in \mathbb{A}(X)$\textit{\ with }$\succ $ $=$ $\vartriangleright $\textit{,
we have }$d(\succsim ,\trianglerighteq )=0$\textit{.}

\bigskip

In the vast majority of the literature on distance functions on preference
relations, it is assumed that the preference relations under consideration
are partial orders. In that setup, or more generally if we wish to define a
metric on the set of all antisymmetric acyclic orders on $X,$ Axiom 4 is
vacuously satisfied.

\subsection{Characterization Theorems}

Let $\succsim $ and $\trianglerighteq $ be two acyclic orders on $X.$ By
Theorem 2.2, we may determine a chain of one-step perturbations that
transform $\succsim $ into $\trianglerighteq $, while Axiom 1 allows us to
find the distance between $\succsim $ and $\trianglerighteq $ by summing up
the distances between each consecutive perturbations in this chain. In turn,
Axiom 3 allows us to compute these distances in terms of rather special
partial orders (of the form $\succsim _{ab}$ and $\succsim _{ab}^{+}$). In
addition, we can compute these distances exactly by using Axioms 1 and 3
jointly.

While there are some technicalities to sort out, this strategy leads to the
following characterization theorem:

\bigskip

\noindent \textbf{Theorem 3.1.} \textit{For any nonempty finite set }$X$%
\textit{, a semimetric }$d:\mathbb{A}(X)\times \mathbb{A}(X)\rightarrow
\lbrack 0,\infty )$\textit{\ satisfies Axioms 1, 3 and 4 if, and only if, }$d
$ \textit{is the }$\mu $\textit{-top-difference semimetric for some positive
measure }$\mu $ \textit{on }$2^{X}$\textit{.}

\bigskip

Adding Axiom 2 to the mix yields:

\bigskip

\noindent \textbf{Theorem 3.2.} \textit{For any nonempty finite set }$X$%
\textit{, a semimetric }$d:\mathbb{A}(X)\times \mathbb{A}(X)\rightarrow
\lbrack 0,\infty )$\textit{\ satisfies Axioms 1-4 if, and only if, }$d$ 
\textit{is the top-difference semimetric.}

\bigskip

The remaining part of this section is devoted to proving Theorem 3.1; the
proof of Theorem 3.2 will be contained within that of Theorem 3.1. To prove
the \textquotedblleft if\textquotedblright\ part of this theorem, we will
use the following fact:

\bigskip

\noindent \textbf{Lemma 3.3.} \textit{For any }$\succsim $, $\succsim _{0}$, 
$\trianglerighteq $\textit{\ }$\in \mathbb{A(}X)$ \textit{with }$\succsim $%
\thinspace $\rightarrow $\thinspace $\succsim _{0}$\thinspace $%
\twoheadrightarrow $\thinspace $\trianglerighteq ,$\textit{\ and }$%
S\subseteq X,$ \textit{the sets }$M(S,\succsim )\triangle M(S,\succsim _{0})$
\textit{and} $M(S,\succsim _{0})\triangle M(S,\trianglerighteq )$ \textit{%
are disjoint, and their union equals }$M(S,\succsim )\triangle
M(S,\trianglerighteq )$.

\begin{proof}
There are two cases to consider. In the first case, there exist $a,b\in X$
such that $\succsim _{0}$ $=$ $\succsim $ $\oplus (a,b),$ $(a,b)\in $ Inc$%
(\succ )$ and $a\vartriangleright b$. In this case, by definition of $%
\succsim _{0},$ we have $a\succ _{0}b.$ Note that if either $a\notin S$ or $%
b\notin M(S,\succsim ),$ we have $M(S,\succsim )=M(S,\succsim _{0}),$ so
there is nothing to prove. Let us then assume that $a\in S$ and $b\in
M(S,\succsim ).$ Since $a\succ _{0}b$ and $a\vartriangleright b$, we then
have $M(S,\succsim )=M(S,\succsim _{0})\sqcup \{b\}$ while $b$ belongs to
neither $M(S,\succsim _{0})$ nor $M(S,\trianglerighteq ).$ It follows that $%
M(S,\succsim )\triangle M(S,\succsim _{0})=\{b\}$ while $b\in M(S,\succsim
)\triangle M(S,\trianglerighteq ).$ But then%
\begin{eqnarray*}
M(S,\succsim _{0})\triangle M(S,\trianglerighteq ) &=&(M(S,\succsim
)\backslash \{b\})\triangle M(S,\trianglerighteq ) \\
&=&(M(S,\succsim )\triangle M(S,\trianglerighteq )\backslash \{b\}\text{.}
\end{eqnarray*}%
The two assertions of the lemma follow from these calculations.

In the second case, there exist $a,b\in X$ such that $\succsim _{0}$ $=$ $%
\succsim $ $\ominus (a,b),$ $a\succ b$, not $a\vartriangleright b$ and (\ref%
{proper}) holds. If either $a\notin S$ or $b\notin M(S,\succsim _{0}),$ we
have $M(S,\succsim )=M(S,\succsim _{0}),$ so there is nothing to prove. We
thus assume $a\in S$ and $b\in M(S,\succsim _{0}).$ Then, since $a\succ b,$ $%
b$ does not belong to $M(S,\succsim ),$ and it readily follows from the
definition of $\succsim _{0}$ that $M(S,\succsim _{0})=M(S,\succsim )\sqcup
\{b\}.$ On the other hand, we now have $b\in M(S,\trianglerighteq )$.
(Otherwise, there exists an $x\in S$ with $x\vartriangleright b,$ so (\ref%
{proper}) implies $x\succ b.$ Given that $a\vartriangleright b$ is not true, 
$x$ must be distinct from $a,$ so we must conclude that $b$ is not $\succsim
_{0}$-maximal in $S,$ a contradiction.) This implies $b\in M(S,\succsim
)\triangle M(S,\trianglerighteq ),$ and therefore, 
\begin{eqnarray*}
M(S,\succsim _{0})\triangle M(S,\trianglerighteq ) &=&(M(S,\succsim )\sqcup
\{b\})\triangle M(S,\trianglerighteq ) \\
&=&(M(S,\succsim )\triangle M(S,\trianglerighteq )\backslash \{b\}\text{.}
\end{eqnarray*}%
The two assertions of the lemma follow from these calculations.
\end{proof}

Let $\mu $ be any positive measure in $2^{X}$. An obvious application of
Lemma 3.3 shows that $D^{\mu }$ satisfies Axiom 1', and by induction, Axiom
1. On the other hand, for any distinct $a,b\in X,$ we have $M(S,\succsim
_{ab})=M(S,\succsim _{ab}^{+})$ for every $S\subseteq X$ distinct from $%
\{a,b\},$ while $M(\{a,b\},\succsim _{ab})=\{a,b\}$ and $M(\{a,b\},\succsim
_{ab}^{+})=\{a\},$ so we obviously have 
\begin{equation}
D^{\mu }(\succsim _{ab},\succsim _{ab}^{+})=\mu \{b\}.  \label{ni}
\end{equation}%
This shows that $D$ satisfies Axiom 2. In turn, to show that $D^{\mu }$
satisfies Axiom 3, take any $\succsim $ $\in \mathbb{A}(X)$ and $a,b\in X.$
Assume first that $a$ and $b$ are not $\succ $-comparable, and put $\succsim
_{0}$ $=$ $\succsim $ $\oplus (a,b).$ As we have shown in the proof of Lemma
3.3, $M(S,\succsim )\triangle M(S,\succsim _{0})=\varnothing $ if either $%
a\notin S$ or $b\notin M(S,\succsim _{0}),$ while $M(S,\succsim )\triangle
M(S,\succsim _{0})=\{b\}$ if $a\in S$ and $b\in M(S,\succsim _{0}).$
Therefore, where $\mathcal{S}:=\{S\in 2^{X}:a\in S$ and $b\in M(S,\succsim
_{0})\},$ we have 
\begin{equation}
D^{\mu }(\succsim ,\succsim _{0})=\sum_{S\in \mathcal{S}}\mu
(\{b\})=\left\vert \mathcal{S}\right\vert \mu (\{b\}).  \label{nic}
\end{equation}%
But, since $a\succ b$ is false, we have $\left\vert \mathcal{S}\right\vert
=2^{N(b,\succsim )-1},$ and combining this with (\ref{ni}) and (\ref{nic}),
we find $D^{\mu }(\succsim ,\succsim _{0})=2^{N(b,\succsim )-1}\mu
\{b\}=2^{N(b,\succsim )-1}D^{\mu }(\succsim _{ab},\succsim _{ab}^{+}),$ as
desired. That $D^{\mu }(\succsim ,\succsim \ominus (a,b))=2^{N(b,\succsim
)}D^{\mu }(\succsim _{ab},\succsim _{ab}^{+})$ when $a\succ b$ is
analogously proved. Finally, it is plain that $D^{\mu }$ satisfies Axiom 4.
We conclude that $D^{\mu }$ satisfies Axioms 1-4.

We now proceed to prove the \textquotedblleft only if\textquotedblright\
part of Theorem 3.1. First, a preliminary observation:

\bigskip

\noindent \textbf{Lemma 3.4.} \textit{Let }$d:\mathbb{A}(X)\times \mathbb{A}%
(X)\rightarrow \lbrack 0,\infty )$\ \textit{be a semimetric that satisfies
Axioms 1 and 3. Then,} 
\begin{equation*}
d(\succsim _{ab},\succsim _{ab}^{+})=d(\succsim _{cb},\succsim _{cb}^{+})%
\hspace{0.2in}\text{\textit{for every distinct} }a,b,c\in X.
\end{equation*}

\begin{proof}
Take any distinct $a,b,c\in X,$ put $Y:=X\backslash \{a,b,c\},$ and consider
the partial orders $\succsim $ and $\trianglerighteq $ on $X$ whose
asymmetric parts are given as 
\begin{equation*}
Y\succ \{a,b,c\}\hspace{0.2in}\text{and\hspace{0.2in}}Y\vartriangleright
\{a,c\}\vartriangleright b.
\end{equation*}%
(In particular, no two distinct element of $Y$ (if any) are comparable by
either $\succsim $ or $\trianglerighteq $.) Then, $\succsim $\thinspace $%
\rightarrow $\thinspace $\succsim \oplus (a,b)$\thinspace $%
\twoheadrightarrow $\thinspace $\trianglerighteq $ so that $d(\succsim
,\trianglerighteq )=d(\succsim ,\succsim \oplus (a,b))+d(\succsim \oplus
(a,b),\trianglerighteq )$ by Axiom 1'. Now by Axiom 3, $d(\succsim ,\succsim
\oplus (a,b))=(2^{2-1})d(\succsim _{ab},\succsim _{ab}^{+}).$ On the other
hand, we have%
\begin{equation*}
\trianglerighteq \text{ }=\text{ }(\succsim \oplus (a,b))\oplus (c,b),
\end{equation*}%
so applying Axiom 3 again yields $d(\succsim \oplus (a,b),\trianglerighteq
)=(2^{1-1})d(\succsim _{cb},\succsim _{cb}^{+}).$ Conclusion:%
\begin{equation*}
d(\succsim ,\trianglerighteq )=2d(\succsim _{ab},\succsim
_{ab}^{+})+d(\succsim _{cb},\succsim _{cb}^{+})\text{.}
\end{equation*}%
But we also have $\succsim $\thinspace $\rightarrow $\thinspace $\succsim
\oplus (c,b)$\thinspace $\twoheadrightarrow $\thinspace $\trianglerighteq $
and $\trianglerighteq $ $=$ $(\succsim \oplus (c,b))\oplus (a,b),$ so
repeating this reasoning yields%
\begin{equation*}
d(\succsim ,\trianglerighteq )=d(\succsim _{ab},\succsim
_{ab}^{+})+2d(\succsim _{cb},\succsim _{cb}^{+})\text{.}
\end{equation*}%
Combining these two equations gives $d(\succsim _{ab},\succsim
_{ab}^{+})=d(\succsim _{cb},\succsim _{cb}^{+})$.
\end{proof}

Now let $d$ be a semimetric on $\mathbb{A}(X)$ that satisfies Axioms 1, 3
and 4. For any $b\in X,$ we define $m_{b}:=d(\succsim _{ab},\succsim
_{ab}^{+})$ where $a\in X\backslash \{b\}.$ By Lemma 3.4, $m_{b}$ is
well-defined nonnegative real number for each $b\in X.$ We define $\mu
:2^{X}\rightarrow \lbrack 0,\infty )$ by $\mu (\varnothing ):=0$ and $\mu
(S):=\sum_{b\in S}m_{b}$ for every nonempty $S\subseteq X.$ Obviously, $\mu $
is a positive measure on $2^{X}$ (and it is the counting measure if $d$
satisfies Axiom 2.) We will complete our proof by showing that $d=D^{\mu }$.

Take any $\succsim $ $\in \mathbb{A}(X).$ Then, for any $(a,b)\in $ Inc$%
(\succ ),$%
\begin{eqnarray*}
d(\succsim ,\succsim \oplus (a,b)) &=&2^{N(b,\succsim )-1}d(\succsim
_{ab},\succsim _{ab}^{+}) \\
&=&2^{N(b,\succsim )-1}\mu (\{b\}) \\
&=&D^{\mu }(\succsim ,\succsim \oplus (a,b)),
\end{eqnarray*}%
where the first equality follows from Axiom 3, the second follows from the
fact that $\mu (\{b\})=m_{b}=d(\succsim _{ab},\succsim _{ab}^{+})$ for any $%
a\in X\backslash \{b\},$ and the third was established above at the end of
the proof of the \textquotedblleft if\textquotedblright\ part of the
theorem. If $a\succ b,$ the analogous reasoning would show instead that $%
d(\succsim ,\succsim \ominus (a,b))=D^{\mu }(\succsim ,\succsim \ominus
(a,b)).$ Conclusion: $d$ and $D^{\mu }$ have the same value at $(\succsim
,\trianglerighteq )$ for every $\succsim ,$ $\trianglerighteq $ $\in \mathbb{%
A}(X)$ where $\trianglerighteq $ is a one-step perturbation of $\succsim $.

Now take any $\succsim ,$ $\trianglerighteq $ $\in \mathbb{A}(X)$ and assume
that the symmetric parts of these relations are the same. If $%
\trianglerighteq $ is a one-step perturbation of $\succsim ,$ we know that $%
d(\succsim ,\trianglerighteq )=D^{\mu }(\succsim ,\trianglerighteq ).$
Otherwise, we apply Theorem 2.2 to find an integer $n\geq 2$ and $\succsim
_{0},...,\succsim _{n-2}\in \mathbb{A}(X)$ such that $\succsim $\thinspace $%
\rightarrow $\thinspace $\succsim _{0}$\thinspace $\twoheadrightarrow $%
\thinspace $\trianglerighteq $ and $\succsim _{k-1}$\thinspace $\rightarrow $%
\thinspace $\succsim _{k}$\thinspace $\twoheadrightarrow $\thinspace $%
\trianglerighteq $ for each $k=1,...,n-1,$ and $\succsim _{n-1}$\thinspace $%
= $\thinspace $\trianglerighteq $. Consequently, applying Axiom 1'
inductively,%
\begin{eqnarray*}
d(\succsim ,\trianglerighteq ) &=&d(\succsim ,\succsim _{0})+\cdot \cdot
\cdot +d(\succsim _{n-2},\succsim _{n-1}) \\
&=&D^{\mu }(\succsim ,\succsim _{0})+\cdot \cdot \cdot +D^{\mu }(\succsim
_{n-2},\trianglerighteq ) \\
&=&D^{\mu }(\succsim ,\trianglerighteq )
\end{eqnarray*}%
where the third equality follows from the fact that $D^{\mu }$ satisfies
Axiom 1'.

Finally, take any $\succsim ,$ $\trianglerighteq $ $\in \mathbb{A}(X)$, and
define $\succsim ^{\ast }$ $:=$ $\succ $ $\sqcup \triangle _{X}$ and $%
\trianglerighteq ^{\ast }$ $:=$ $\vartriangleright $ $\sqcup \triangle _{X}$%
. Then, $\succsim ^{\ast },\trianglerighteq ^{\ast }\in \mathbb{A}(X)$ and $%
d(\succsim ^{\ast },\trianglerighteq ^{\ast })=D^{\mu }(\succsim ^{\ast
},\trianglerighteq ^{\ast })$ by what we have found in the previous
paragraph. But, by Axiom 4, $d(\succsim ,\succsim ^{\ast
})=0=d(\trianglerighteq ,\trianglerighteq ^{\ast })$. Since $d$ is a
semimetric, therefore, 
\begin{equation*}
d(\succsim ,\trianglerighteq )=d(\succsim ,\succsim ^{\ast })+d(\succsim
^{\ast },\trianglerighteq ^{\ast })+d(\trianglerighteq ^{\ast
},\trianglerighteq )=d(\succsim ^{\ast },\trianglerighteq ^{\ast })=D^{\mu
}(\succsim ^{\ast },\trianglerighteq ^{\ast }).
\end{equation*}%
Since $M(S,\succsim ^{\ast })=M(S,\succsim )$ and $M(S,\trianglerighteq
^{\ast })=M(S,\trianglerighteq )$ for every $S\subseteq X,$ we have $D^{\mu
}(\succsim ^{\ast },\trianglerighteq ^{\ast })=D^{\mu }(\succsim
,\trianglerighteq ),$ and hence obtains $d(\succsim ,\trianglerighteq
)=D^{\mu }(\succsim ,\trianglerighteq ).$ The proof of Theorem 3.1 is now
complete.

\subsection{Top-Difference Metrics vs. Weighted KSB Metrics}

As we have noted in Section 1, Can \cite{Can} and Hassanzadeh and Milenkovic 
\cite{Has-M} were motivated by observations such as the one we presented in
Example 1.1, and have consequently proposed a class of metrics that consist
of weighted forms of the classical Kemeny-Snell metric. It should be noted
that these metrics are defined only on $\mathbb{L}(X)$, the set of all 
\textit{linear} orders on $X.$ Moreover, it is not at all clear how to
extend these metrics (axiomatically or even simply by definition) to the
domains like $\mathbb{P}(X)$ or $\mathbb{P}^{\ast }(X)$. As such, we can
make a comparison with these metrics and the $\mu $-top-difference
semimetrics only by restricting the domain of the latter to $\mathbb{L}(X)$.
(As noted earlier, on this domain, any $D^{\mu }$ acts as a metric.)

Let $n:=\left\vert X\right\vert ,$ and let $\Sigma $ denote the set of all
permutations $\sigma $ on $\{1,...,n\}$ for which there is a $k\in
\{1,...,n-1\}$ with $\sigma (k)=k+1,$ $\sigma (k+1)=\sigma (k),$ and $\sigma
(i)=i$ for all $i\neq k,k+1.$ (\cite{Has-M} refer to such permutations as 
\textit{adjacent transpositions.}) In what follows, we abuse notation and
write $\sigma (x_{1},...,x_{n})$ for the $n$-vector $(x_{\sigma
(1)},...,x_{\sigma (n)})$ for any $x_{1},...,x_{n}\in X$ and $\sigma \in
\Sigma $.

Next, for any $\succsim $ $\in \mathbb{L}(X),$ let us agree to write $%
v(\succsim )$ for the $n$-vector $(x_{1},...,x_{n})$ where $x_{1}\succ \cdot
\cdot \cdot \succ x_{n}$. Finally, for any (weight function) $\omega :\Sigma
\rightarrow \lbrack 0,\infty ),$ we define the real map $d_{\omega }$ on $%
\mathbb{L}(X)\times \mathbb{L}(X)$ by%
\begin{equation*}
d_{\omega }(\succsim ,\trianglerighteq ):=\min \sum_{i=1}^{k}\omega (\sigma
_{i})
\end{equation*}%
where the minimum is taken over all $k\in \mathbb{N}$ and $\sigma
_{1},...,\sigma _{k}\in \Sigma $ such that $v(\trianglerighteq )=(\sigma
_{1}\circ \cdot \cdot \cdot \circ \sigma _{k})v(\succsim )$. It is easy to
check that this is indeed a metric on $\mathbb{L}(X);$ it is referred to as
a \textit{weighted Kendall metric} by \cite{Has-M}. Of course, for the
weight function $\omega $ that equals 2 everywhere, $d_{\omega }$ becomes
precisely the Kemeny-Snell-Bogart metric on $\mathbb{L}(X)$.

The first observation we would like to make here is that if we choose the
weight function $\omega $ such that $\omega (\sigma )=2^{n-k}$ for the
adjacent transposition $\sigma $ with $\sigma (k)=k+1$ and $\sigma
(k+1)=\sigma (k),$ then $d_{\omega }$ becomes identical to the
top-difference metric on $\mathbb{L}(X)$, that is, $d_{\omega }=D|_{\mathbb{L%
}(X)\times \mathbb{L}(X)}$ for this special weight function $\omega .$ (We
omit the straightforward proof.) Second, when $n\geq 3,$ there is no weight
function $\omega $ such that $d_{\omega }=D^{\mu }|_{\mathbb{L}(X)\times 
\mathbb{L}(X)}$ unless $\mu $ is indeed the counting measure. To see this,
pick any measure $\mu $ on $2^{X}$ such that $\mu (\{x\})\neq \mu (\{y\})$
for some $x,y\in X.$ Now take any $a\in X\backslash \{x,y\},$ and consider
the linear orders $\succsim _{1},...,\succsim _{4}$ on $X$ such that%
\begin{equation*}
x\succ _{1}\cdot \cdot \cdot \succ _{1}a\succ _{1}y\text{\hspace{0.1in}and%
\hspace{0.1in}}x\succ _{2}\cdot \cdot \cdot \succ _{2}y\succ _{2}a\text{,}
\end{equation*}%
and%
\begin{equation*}
y\succ _{3}\cdot \cdot \cdot \succ _{3}x\succ _{3}a\text{,\hspace{0.1in}and }%
y\succ _{4}\cdot \cdot \cdot \succ _{4}a\succ _{4}x,
\end{equation*}%
with the understanding that the unspecified parts of all of these linear
orders agree. Then, it is plain that $d_{\omega }(\succsim _{1},\succsim
_{2})=d_{\omega }(\succsim _{3},\succsim _{4})$ for any $\omega :\Sigma
\rightarrow \lbrack 0,\infty ).$ By contrast, $D^{\mu }(\succsim
_{1},\succsim _{2})=\mu (\{a,y\})\neq \mu (\{a,x\})=D^{\mu }(\succsim
_{3},\succsim _{4})$, which shows that $D^{\mu }$ is distinct from $%
d_{\omega }$ no matter how we may choose the weight function $\omega $. We
proved:

\bigskip

\noindent \textbf{Proposition 3.5.} \textit{For any finite set }$X$ \textit{%
with} $\left\vert X\right\vert \geq 3$\textit{, the only }$\mu $\textit{%
-top-difference metric which is also a weighted Kendall metric is the
top-difference metric on }$\mathbb{L}(X).$

\bigskip 

Thus, not only does the metrization approach we develop here applies well
beyond $\mathbb{L}(X),$ even on this space it is quite distinct from those
of \cite{Can} and \cite{Has-M}. There is only one exception to this, namely,
the top-difference metric $D$ on $\mathbb{L}(X).$ This metric is the only
one that lies in the intersection of the $D^{\mu }$ class and the class of
metrics introduced by \cite{Can} and \cite{Has-M}.\footnote{%
We should note that the special weight function that yields $D$ within the
class of all weighted Kendall metrics is not mentioned in either \cite{Can}
or \cite{Has-M}. These papers do not provide an axiomatization for this
semimetric even on $\mathbb{L}(X).$} This observation further singles out
this semimetric as the most important member of the class of $D^{\mu }$
metrics, and provides motivation for its further investigation.

\subsection{On the Computational Complexity for $D^{\protect\mu }$}

While the intuition behind the $D^{\mu }$ metrics appears convincing, and it
is reinforced by Theorems 3.1 and 3.2, the computation of distances between
two acyclic preferences $\succsim $ and $\trianglerighteq $ on $X$ according
to any one of these metrics require one compute the symmetric difference
between $M(S,\succsim )$ and $M(S,\trianglerighteq )$ for all subsets $S$ of 
$X.$ As this set is empty whenever $\left\vert S\right\vert \leq 1,$ this
means we have to compute $M(S,\succsim )\triangle M(S,\trianglerighteq )$
for $2^{\left\vert X\right\vert }-\left\vert X\right\vert -1$ many subsets
of $X.$ As $\left\vert X\right\vert $ gets larger, this becomes a
computationally daunting task. This is in stark contrast with the
computation of distances relative to the Kemeny-Snell-Bogart metric which
requires at most polynomial time with respect to the size of $X.$

Fortunately, there is a more efficient way of computing $D^{\mu }(\succsim
,\trianglerighteq )$ for any given $\succsim ,\trianglerighteq $ $\in 
\mathbb{A}(X)$ which we now explore. For any $S\subseteq X,$ let us first
write%
\begin{equation*}
\triangle _{S}(\succsim ,\trianglerighteq ):=M(S,\succsim )\triangle
M(S,\trianglerighteq )
\end{equation*}%
to simplify our notation. Then, for any fixed positive measure $\mu $ on $%
2^{X},$ we have%
\begin{eqnarray*}
D^{\mu }(\succsim ,\trianglerighteq ) &=&\sum_{S\subseteq X}\mu (\triangle
_{S}(\succsim ,\trianglerighteq )) \\
&=&\sum_{S\subseteq X}\sum_{x\in S}\mu (\{x\})\mathbf{1}_{\triangle
_{S}(\succsim ,\trianglerighteq )}(x) \\
&=&\sum_{x\in X}\left( \sum_{\substack{ S\subseteq X \\ S\ni x}}\mathbf{1}%
_{\triangle _{S}(\succsim ,\trianglerighteq )}(x)\right) \mu (\{x\})
\end{eqnarray*}%
In other words,%
\begin{equation*}
D^{\mu }(\succsim ,\trianglerighteq )=\sum_{x\in X}\theta _{x}(\succsim
,\trianglerighteq )\mu (\{x\})
\end{equation*}%
where $\theta _{x}(\succsim ,\trianglerighteq )$ is the number of all
subsets $S$ of $X$ such that $x\in \triangle _{S}(\succsim ,\trianglerighteq
)$.

Let us now fix any $x\in X,$ and calculate $\theta _{x}(\succsim
,\trianglerighteq )$. To this end, let us define the following three sets:%
\begin{equation*}
A_{x}(\succsim ,\trianglerighteq ):=\{a\in X\backslash \{x\}:\text{ not }%
a\succ x\text{ and not }a\vartriangleright x\},
\end{equation*}%
and 
\begin{equation*}
B_{x}(\succsim ,\trianglerighteq ):=\{a\in X\backslash \{x\}:a\succ x\text{
but not }a\vartriangleright x\},
\end{equation*}%
and%
\begin{equation*}
C_{x}(\succsim ,\trianglerighteq ):=\{a\in X\backslash
\{x\}:a\vartriangleright x\text{ but not }a\succ x\}\text{.}
\end{equation*}%
We denote the cardinality of the first of these sets by $\alpha
_{x}(\succsim ,\trianglerighteq )$. Notice first that $x\in M(S,\succsim
)\backslash M(S,\trianglerighteq )$ iff $S=\{x\}\sqcup K\sqcup L$ for some $%
K\subseteq A_{x}(\succsim ,\trianglerighteq )$ and some \textit{nonempty} $%
L\subseteq C_{x}(\succsim ,\trianglerighteq ).$ There are exactly $2^{\alpha
_{x}(\succsim ,\trianglerighteq )}(2^{\left\vert C_{x}(\succsim
,\trianglerighteq )\right\vert }-1)$ many such sets. On the other hand, by
the same logic, there are $2^{\alpha _{x}(\succsim ,\trianglerighteq
)}(2^{\left\vert B_{x}(\succsim ,\trianglerighteq )\right\vert }-1)$ many
subsets $S$ of $X$ such that $x\in M(S,\trianglerighteq )\backslash
M(S,\succsim ).$ It follows that%
\begin{equation*}
\theta _{x}(\succsim ,\trianglerighteq )=2^{\alpha _{x}(\succsim
,\trianglerighteq )}(2^{\left\vert B_{x}(\succsim ,\trianglerighteq
)\right\vert }+2^{\left\vert C_{x}(\succsim ,\trianglerighteq )\right\vert
}-2).
\end{equation*}%
Next, notice that $A_{x}(\succsim ,\trianglerighteq )\sqcup B_{x}(\succsim
,\trianglerighteq )=\{a\in X\backslash \{x\}:$ not $a\vartriangleright x\},$
whence 
\begin{equation*}
\alpha _{x}(\succsim ,\trianglerighteq )+\left\vert B_{x}(\succsim
,\trianglerighteq )\right\vert =n-\left\vert x^{\uparrow ,\vartriangleright
}\right\vert -1
\end{equation*}%
where $n:=\left\vert X\right\vert $, and as we defined in Section 2.1, $%
x^{\uparrow ,\vartriangleright }$ is the principal ideal of $x$ with respect
to $\vartriangleright $. Of course, the analogous reasoning shows that $%
\alpha _{x}(\succsim ,\trianglerighteq )+\left\vert C_{x}(\succsim
,\trianglerighteq )\right\vert =n-\left\vert x^{\uparrow ,\succ }\right\vert
-1$ as well. Consequently,%
\begin{equation*}
\theta _{x}(\succsim ,\trianglerighteq )=2^{n-\left\vert x^{\uparrow
,\vartriangleright }\right\vert -1}+2^{n-\left\vert x^{\uparrow ,\succ
}\right\vert -1}-2^{\alpha _{x}(\succsim ,\trianglerighteq )+1}.
\end{equation*}

Combining the computations of the previous two paragraphs yields an
alternative method of calculating the distance between $\succsim $ and $%
\trianglerighteq $ with respect to $D^{\mu }$:%
\begin{equation}
D^{\mu }(\succsim ,\trianglerighteq )=\sum_{x\in X}\left[ 2^{n-\left\vert
x^{\uparrow ,\vartriangleright }\right\vert -1}+2^{n-\left\vert x^{\uparrow
,\succ }\right\vert -1}-2^{\alpha _{x}(\succsim ,\trianglerighteq )+1}\right]
\mu (\{x\}).  \label{new}
\end{equation}%
This formula does not look particularly appealing at first glance. It is not
even clear that it defines a semimetric on $\mathbb{A}(X),$ and it is
certainly not intuitive. However, it has a significant computational
advantage over the formula we defined $D^{\mu }$ with. Indeed, this formula
uses only \textquotedblleft local\textquotedblright\ knowledge about the
involved acyclic orders. As a consequence, the computation of the numbers $%
\left\vert x^{\uparrow ,\vartriangleright }\right\vert ,$ $\left\vert
x^{\uparrow ,\succ }\right\vert $ and $\alpha _{x}(\succsim
,\trianglerighteq )$ for each $x\in X,$ and hence the above formula, take at
most polynomial time with respect to the size of $X,$ which parallels the
computational efficiency of the Kemeny-Snell-Bogart metric. Any sort of a
program that is primed to compute the values of $D^{\mu }$ should thus
utilize (\ref{new}) instead of (\ref{old}). The computational superiority of
(\ref{new}) over (\ref{old}) will be further witnessed in the next
subsection.

\subsection{The Distance Between Linear Orders}

The family $\mathbb{L}(X)$ of linear orders on $X$ arises in numerous
applications, ranging from voting theory to stable matching, random utility
theory, etc.. Indeed, the Kemeny-Snell-Bogart metric is primarily applied on 
$\mathbb{L}(X)$ (and as such, it is often simply referred to as the \textit{%
Kemeny-Snell metric}). It is thus natural to ask if there is an easy way of
computing the top-difference metric $D$ on $\mathbb{L}(X)\times \mathbb{L}(X)
$. (We recall that $D$ is a metric on $\mathbb{L}(X),$ not only a
semimetric.) We next provide such a formula by using (\ref{new}).

Take any $\succsim ,\trianglerighteq $ $\in \mathbb{L}(X)$, and put $%
n:=\left\vert X\right\vert $. Given that $\succsim $ is a linear order, for
every $i\in \{0,...,n-1\},$ there is a unique $x\in X$ such that $\left\vert
x^{\downarrow ,\succ }\right\vert =i.$ Moreover, again by linearity, $%
\left\vert x^{\downarrow ,\succ }\right\vert =n-\left\vert x^{\uparrow
,\succ }\right\vert -1$ for each $x\in X.$ It follows that 
\begin{equation*}
\sum_{x\in X}2^{n-\left\vert x^{\uparrow ,\succ }\right\vert
-1}=\sum_{i=0}^{n-1}2^{i}=2^{n}-1.
\end{equation*}%
Since, analogously, we also have $\sum_{x\in X}2^{n-\left\vert x^{\uparrow
,\vartriangleright }\right\vert -1}=2^{n}-1,$ the formula (\ref{new}) yields%
\begin{equation*}
D(\succsim ,\trianglerighteq )=2(2^{n}-1)-\sum_{x\in X}2^{\alpha
_{x}(\succsim ,\trianglerighteq )+1}\text{.}
\end{equation*}%
Next, notice that $\alpha _{x}(\succsim ,\trianglerighteq )$ is none other
than the number of all elements of $X$ that are strictly below $x$ with
respect to both $\succsim $ and $\trianglerighteq $ (again because $\succsim 
$ and $\trianglerighteq $ are linear orders on $X$). Consequently, we arrive
at%
\begin{equation*}
D(\succsim ,\trianglerighteq )=2(2^{n}-1)-\sum_{x\in X}2^{\left\vert
x^{\downarrow ,\succ }\cap x^{\downarrow ,\vartriangleright }\right\vert +1}%
\text{.}
\end{equation*}%
This shows that to find the distance between two linear orders on $X,$ all
one has to do is to count the elements in the intersection of the principal
filters of each $x\in X$ with respect to $\succsim $ and $\trianglerighteq $%
. This is very efficient, as it allows us to work with the orders $\succsim $
and $\trianglerighteq $ separately.

\section{Diameter of the Preference Space $(\mathbb{A}(X),D)$}

To get a better sense of the \textquotedblleft distance\textquotedblright\
between two preference relations in practice, one should really have a basic
benchmark. In particular, it may be useful to know the \textit{diameter} of
the space of preferences one is interested in with respect to the semimetric
at hand. In this section we thus attempt to get some simple lower estimates
for the diameter of $\mathbb{A}(X)$ and $\mathbb{P}(X)$ with respect to the
top-difference semimetric $D.$ (We denote the diameter operator relative to $%
D$ by diam$_{D}(\cdot )$.)

Let us denote the cardinality of $X$ by $n;$ recall that $n\geq 2$. The
diameter problem is easily treated in the case of linear orders. Indeed, for
any $\succsim ,$ $\succsim ^{\prime }$ $\in \mathbb{L}(X),$ the cardinality
of $M(S,\succsim )\triangle M(S,\succsim ^{\prime })$ is at most 2 for any $%
S\subseteq X$ with at least two elements. Therefore, the largest possible
value for $D(\succsim ,\succsim ^{\prime })$ is $2$ times the number of all $%
S\subseteq X$ with $\left\vert S\right\vert \geq 2,$ namely, $2(2^{n}-n-1).$
But if we enumerate $X$ as $\{x_{1},...,x_{n}\},$ and choose $\succsim $ and 
$\succsim ^{\prime }$ orthogonally to each other as $x_{1}\succ \cdot \cdot
\cdot \succ x_{n}$ and $x_{n}\succ ^{\prime }\cdot \cdot \cdot \succ
^{\prime }x_{1},$ then $\left\vert M(S,\succsim )\triangle M(S,\succsim
^{\prime })\right\vert =2$ for all $S\subseteq X$ with $\left\vert
S\right\vert \geq 2$. Thus:%
\begin{equation}
\text{diam}_{D}(\mathbb{L}(X))=2(2^{n}-n-1).  \label{di1}
\end{equation}%
To put this number in some perspective, we report its value in the table
below for the first nine values of $n,$ next to the cardinality $n!$ of $%
\mathbb{L}(X).$

The situation is more complicated for total preorders. To examine this case,
we fix any $m\in \{1,...,n-1\},$ and consider the total preorders $\succsim $
and $\succsim ^{\prime }$ on $X$ such that%
\begin{equation*}
x_{1}\sim \cdot \cdot \cdot \sim x_{m}\succ x_{m+1}\succ \cdot \cdot \cdot
\succ x_{n}
\end{equation*}%
and%
\begin{equation*}
x_{m+1}\sim ^{\prime }\cdot \cdot \cdot \sim ^{\prime }x_{n}\succ ^{\prime
}x_{1}\succ ^{\prime }\cdot \cdot \cdot \succ ^{\prime }x_{m}\text{.}
\end{equation*}%
Now let $A:=\{x_{1},...,x_{m}\}$ and $B:=\{x_{m+1},...,x_{n}\},$ and note
that%
\begin{equation*}
\left\vert M(S,\succsim )\triangle M(S,\succsim ^{\prime })\right\vert
=\left\{ 
\begin{array}{ll}
\left\vert S\right\vert -1, & \text{if }S\subseteq A\text{ or }S\subseteq B,
\\ 
\left\vert S\right\vert , & \text{otherwise}%
\end{array}%
\right.
\end{equation*}%
for any $S\subseteq X.$ Where $\mathcal{S}:=\{S\in 2^{X}\backslash
\{\varnothing \}:S\cap A\neq \varnothing \neq S\cap B\},$ we thus have%
\begin{eqnarray*}
D(\succsim ,\succsim ^{\prime }) &=&\sum_{S\in \mathcal{S}}\left\vert
S\right\vert +\sum_{\varnothing \neq S\subseteq A}(\left\vert S\right\vert
-1)+\sum_{\varnothing \neq S\subseteq B}(\left\vert S\right\vert -1) \\
&=&\sum_{\varnothing \neq S\subseteq X}\left\vert S\right\vert
-\sum_{\varnothing \neq S\subseteq A}1-\sum_{\varnothing \neq S\subseteq A}1
\\
&=&\sum_{k=1}^{n}k\binom{n}{k}-\left\vert 2^{A}\backslash \{\varnothing
\}\right\vert -\left\vert 2^{B}\backslash \{\varnothing \}\right\vert \\
&=&\sum_{k=1}^{n}k\binom{n}{k}+2-2^{m}-2^{n-m}\text{.}
\end{eqnarray*}%
It is readily checked that $t\mapsto 2^{t}+2^{n-t}$ is a symmetric and
strictly convex function on $[0,n];$ this function attains its unique global
minimum at $\frac{n}{2}$. It follows that the map $m\mapsto 2^{m}+2^{m-t}$
achieves its minimum on $\{0,...,m\}$ at $\lfloor \frac{n}{2}\rfloor $.
Combining this fact with the calculation above, and recalling that $%
\sum_{k=1}^{n}k\binom{n}{k}=n2^{n-1}$ (which is easily verified by induction
on $n$) and $\lceil \frac{n}{2}\rceil =n-\lfloor \frac{n}{2}\rfloor $, we
find that $n2^{n-1}+2-2^{\lfloor \frac{n}{2}\rfloor }-2^{\lceil \frac{n}{2}%
\rceil }$ is a lower bound for diam$_{D}(\mathbb{P}_{\text{total}}(X)).$ In
our next result we prove that this lower bound is actually attained.

\bigskip

\noindent \textbf{Theorem 4.1.} \textit{Let }$X$ \textit{be a finite set
with }$n:=\left\vert X\right\vert \geq 2.$ \textit{Then, }%
\begin{equation}
\text{diam}_{D}(\mathbb{P}_{\text{total}}(X))=n2^{n-1}+2-2^{\lfloor \frac{n}{%
2}\rfloor }-2^{\lceil \frac{n}{2}\rceil }\text{.}  \label{di2}
\end{equation}

\noindent \textit{Proof. }Let us begin by noting that for $n=2$ and $n=3,$
it is readily checked that diam$_{D}(\mathbb{L}(X))=$ diam$_{D}(\mathbb{P}_{%
\text{total}}(X))$ and that the right-hand sides of (\ref{di1}) and (\ref%
{di2}) are the same. As this observation readily yields the present theorem
for $n\in \{2,3\},$ we assume $n\geq 4$ in the rest of the proof.

Now define $\eta :\{1,...,n-1\}\rightarrow (-\infty ,0)$ by $\eta
(m):=2-2^{m}-2^{n-m}.$ We have seen above that $\eta (\lfloor \tfrac{n}{2}%
\rfloor )\geq \eta (m)$ for each $m=1,...,n-1,$ and that%
\begin{equation*}
\text{diam}_{D}(\mathbb{P}_{\text{total}}(X))\geq n2^{n-1}+\eta (\lfloor 
\tfrac{n}{2}\rfloor )\text{.}
\end{equation*}%
To prove the converse inequality, we take any total preorders $\succsim $
and $\succsim ^{\prime }$ on $X.$ We must show that $D(\succsim ,\succsim
^{\prime })\leq n2^{n-1}+\eta (\lfloor \tfrac{n}{2}\rfloor )$.

Let us first assume that there is at least one element that is maximal in $X$
with respect to both $\succsim $ and $\succsim ^{\prime }$. Let $\mathcal{A}$
stand for the set of all subsets of $X$ that contain this element, and note
that $\left\vert \mathcal{A}\right\vert =2^{n-1}=-\eta (1).$ Then, $%
\left\vert M(S,\succsim )\triangle M(S,\succsim ^{\prime })\right\vert $ is
at most $\left\vert S\right\vert -1$ for every $S\in \mathcal{A}$ while it
is trivially less than $\left\vert S\right\vert $ for any $S\subseteq X$.
Consequently,%
\begin{eqnarray*}
D(\succsim ,\succsim ^{\prime }) &\leq &\sum_{S\in \mathcal{A}}(\left\vert
S\right\vert -1)+\sum_{S\in 2^{X}\backslash \mathcal{A}}\left\vert
S\right\vert \\
&=&\sum_{S\subseteq X}\left\vert S\right\vert -\left\vert \mathcal{A}%
\right\vert \\
&=&n2^{n-1}+\eta (1) \\
&\leq &n2^{n-1}+\eta (\lfloor \tfrac{n}{2}\rfloor ),
\end{eqnarray*}%
as desired.\footnote{%
The third equality here holds because $\sum_{S\subseteq X}\left\vert
S\right\vert =\sum_{k=1}^{n}k\binom{n}{k}=n2^{n-1}$.}

It remains to consider the case $M(X,\succsim )\cap M(X,\succsim ^{\prime
})=\varnothing $. There are two possibilities to consider in this case.
First, assume that $M(X,\succsim )\sqcup M(X,\succsim ^{\prime })=X.$ In
this case, we put $m:=\left\vert M(X,\succsim )\right\vert ,$ and note that $%
\left\vert M(X,\succsim ^{\prime })\right\vert =n-m$. Let $\mathcal{A}$
stand for the set of all nonempty subsets $S$ of $X$ such that either $%
S\subseteq M(X,\succsim )$ or $S\subseteq \left\vert M(X,\succsim ^{\prime
})\right\vert $. Since $M(X,\succsim )$ and $M(X,\succsim ^{\prime })$ are
disjoint, we have $\left\vert \mathcal{A}\right\vert
=(2^{m}-1)+(2^{n-m}-1)=-\eta (m).$ On the other hand, again, $\left\vert
M(S,\succsim )\triangle M(S,\succsim ^{\prime })\right\vert \leq \left\vert
S\right\vert -1$ for every $S\in \mathcal{A}$. Therefore, carrying out the
same calculation we have done in the previous paragraph yields $D(\succsim
,\succsim ^{\prime })\leq n2^{n-1}+\eta (m)\leq n2^{n-1}+\eta (\lfloor 
\tfrac{n}{2}\rfloor )$, as desired.

The only remaining case is where $M(X,\succsim )\cap M(X,\succsim ^{\prime
})=\varnothing $ and $M(X,\succsim )\sqcup M(X,\succsim ^{\prime })\neq X.$
In this case, to simplify our notation, we put $A:=M(X,\succsim )$, $%
B:=M(X,\succsim ^{\prime })$ and $C=X\backslash (A\sqcup B)$. Let $%
m_{1}:=\left\vert A\right\vert ,$ $m_{2}:=\left\vert B\right\vert ,$ and
note that $\left\vert C\right\vert =n-m_{1}-m_{2}>0.$ Next, we define $%
\mathcal{A}$ exactly as in the previous paragraph, and note that $\left\vert 
\mathcal{A}\right\vert =(2^{m_{1}}-1)+(2^{m_{2}}-1)$ and $\left\vert
M(S,\succsim )\triangle M(S,\succsim ^{\prime })\right\vert \leq \left\vert
S\right\vert -1$ for every $S\in \mathcal{A}$. Finally, we define%
\begin{equation*}
\mathcal{B}:=\{S\in 2^{X}:S\cap A\neq \varnothing ,\text{ }S\cap B\neq
\varnothing \text{ and }S\cap C\neq \varnothing \}\text{.}
\end{equation*}%
Then, 
\begin{subequations}
\begin{eqnarray*}
D(\succsim ,\succsim ^{\prime }) &\leq &\sum_{S\in \mathcal{A}}(\left\vert
S\right\vert -1)+\sum_{S\in \mathcal{B}}(\left\vert S\right\vert -\left\vert
S\cap C\right\vert )+\sum_{S\in 2^{X}\backslash (\mathcal{A}\sqcup \mathcal{B%
})}\left\vert S\right\vert \\
&=&\sum_{S\subseteq X}\left\vert S\right\vert -\left\vert \mathcal{A}%
\right\vert -\sum_{S\in \mathcal{B}}\left\vert S\cap C\right\vert \\
&=&n2^{n-1}-(2^{m_{1}}-1)-(2^{m_{2}}-1)-\sum_{S\in \mathcal{B}}\left\vert
S\cap C\right\vert \text{.}
\end{eqnarray*}%
On the other hand, by definition of $\mathcal{B}$, 
\end{subequations}
\begin{eqnarray*}
\sum_{S\in \mathcal{B}}\left\vert S\cap C\right\vert
&=&(2^{m_{1}}-1)(2^{m_{2}}-1)\sum_{k=1}^{n-m_{1}-m_{2}}k\binom{n-m_{1}-m_{2}%
}{k} \\
&=&(n-m_{1}-m_{2})(2^{m_{1}}-1)(2^{m_{2}}-1)2^{n-m_{1}-m_{2}-1}\text{.}
\end{eqnarray*}%
If $n-m_{1}-m_{2}=1,$ therefore, 
\begin{eqnarray*}
(2^{m_{1}}-1)+(2^{m_{2}}-1)+\sum_{S\in \mathcal{B}}\left\vert S\cap
C\right\vert &=&(2^{m_{1}}-1)+(2^{m_{2}}-1)+(2^{m_{1}}-1)(2^{m_{2}}-1) \\
&=&2^{n-1}-1 \\
&\geq &2^{n-2}+2 \\
&=&-\eta (2).
\end{eqnarray*}%
(The inequality here follows because $n\geq 4$ and $2^{t-1}-2^{t-2}-3\geq 0$
for every $t\geq 4.\footnote{%
This follows from the fact that the map $t\mapsto 2^{t-1}-2^{t-2}-3$ is
(strictly) increasing on $[4,\infty )$ and its value at $4$ is positive.}$)
If, on the other hand, $n-m_{1}-m_{2}\geq 2,$ we have%
\begin{eqnarray*}
\sum_{S\in \mathcal{B}}\left\vert S\cap C\right\vert &\geq
&2(2^{m_{1}}-1)(2^{m_{2}}-1)2^{n-m_{1}-m_{2}-1} \\
&\geq &2^{m_{1}-1}2^{m_{2}-1}2^{n-m_{1}-m_{2}} \\
&=&2^{n-2}\text{.}
\end{eqnarray*}%
(Here we use the fact that $2^{t}-2^{t-1}-1\geq 0$ for every $t\geq 1.$)
Thus, again, we find%
\begin{eqnarray*}
(2^{m_{1}}-1)+(2^{m_{2}}-1)+\sum_{S\in \mathcal{B}}\left\vert S\cap
C\right\vert &\geq &2^{m_{1}}+2^{m_{2}}-2+2^{n-2} \\
&\geq &4-2+2^{n-2} \\
&\geq &2+2^{n-2} \\
&=&-\eta (2).
\end{eqnarray*}%
Returning to the computation of $D(\succsim ,\succsim ^{\prime })$, we then
get 
\begin{subequations}
\begin{eqnarray*}
D(\succsim ,\succsim ^{\prime }) &\leq
&n2^{n-1}-(2^{m_{1}}-1)-(2^{m_{2}}-1)-\sum_{S\in \mathcal{B}}\left\vert
S\cap C\right\vert \\
&\leq &n2^{n-1}+\eta (2) \\
&\leq &n2^{n-1}+\eta (\lfloor \tfrac{n}{2}\rfloor ).
\end{eqnarray*}%
The proof of Theorem 4.1 is now complete. \hspace{0.2in}\hspace{0.2in}%
\hspace{0.2in}\hspace{0.2in}\hspace{0.2in}\hspace{0.2in}\hspace{0.2in}%
\hspace{0.2in}\hspace{0.2in}\hspace{0.2in}\hspace{0.2in} $\square $

\bigskip

For any integer $n\geq 2,$ let us denote the number of total preorders on
the $n$-element set $X$ by $p(n).$ It is a well known combinatorial fact
that this number can be computed as 
\end{subequations}
\begin{equation*}
p(n)=\sum_{k=0}^{n}k!S(n,k)
\end{equation*}%
where $S(n,k)$ is the number of ways an $n$-element set can be partitioned
into $k$ many nonempty sets; these numbers are known as the \textit{Stirling
numbers of the second kind}. Table 1 provides a comparison between $p(n)$
and the $D$-diameter of $\mathbb{P}_{\text{total}}(X)$ up to $n=10$.

\begin{eqnarray*}
&&%
\begin{array}{ccccc}
& \text{diam}_{D}(\mathbb{P}_{\text{total}}(X)) & p(n) & \text{diam}_{D}(%
\mathbb{L}(X)) & n! \\ 
n=2 & 2 & 3 & 2 & 2 \\ 
n=3 & 8 & 13 & 8 & 6 \\ 
n=4 & 26 & 75 & 22 & 24 \\ 
n=5 & 70 & 541 & 52 & 120 \\ 
n=6 & 178 & 4,683 & 114 & 720 \\ 
n=7 & 426 & 47,293 & 240 & 5,040 \\ 
n=8 & 994 & 545,835 & 494 & 40,320 \\ 
n=9 & 2,258 & 7,087,261 & 1,004 & 362,880 \\ 
n=10 & 5,058 & 102,247,563 & 2,026 & 3,628,800%
\end{array}
\\
&& \\
&&\hspace{0.2in}\hspace{0.2in}\hspace{0.2in}\hspace{0.2in}\hspace{0.2in}%
\hspace{0.2in}\hspace{0.2in}\hspace{0.2in}\hspace{0.2in}\text{\textbf{Table 1%
}}
\end{eqnarray*}

\noindent This table suggests that, relative to the size of $\mathbb{P}_{%
\text{total}}(X),$ the $D$-diameter of $\mathbb{P}_{\text{total}}(X)$
remains fairly modest, just as in the case of $\mathbb{L}(X).$

In passing, we note that as an immediate consequence of Theorem 4.1, we have%
\begin{equation*}
\text{diam}_{D}(\mathbb{A}(X))\geq \text{diam}_{D}(\mathbb{P}(X))\geq
n2^{n-1}+2-2^{\lfloor \frac{n}{2}\rfloor }-2^{\lceil \frac{n}{2}\rceil }%
\text{.}
\end{equation*}%
We do not presently know whether or not either of these inequalities hold as
equalities.

\section{On Best Transitive Approximations}

As an acyclic order $\succsim $ on $X$ need not be transitive, a natural
problem is to identify the set of all preorders on $X$ that best
approximates $\succsim $ in the sense of distance minimizing where we
measure distance by $D$ (or by $D^{\mu }$ for some suitable $\mu $). Put
differently, the problem is to compute the metric projection of $\succsim $
in $\mathbb{P}(X)$ relative to $D$ (or $D^{\mu }$). This seems like an
interesting problem, and it should eventually be studied from an algorithmic
perspective. Here we offer a partial solution to it.

First, we simplify the problem by assuming $\succsim $ is antisymmetric.
Second, we concentrate on finding the best approximation to $\succsim $
among all preorders that extend $\succsim $. Recall that a binary relation $%
R $ on $X$ \textit{extends} $\succsim $ if it is reflexive and satisfies $%
\succ $ $\subseteq $ $R^{>}$. (That is, an extension $R$ of $\succsim $ is
particularly faithful to $\succsim $ in that its ranking of any two $%
\succsim $-comparable alternatives is identical to the ranking of those
alternatives by $\succsim $.) We denote the set of all transitive extensions
of $\succsim $ by Ext$(\succsim ).$ For any given positive measure $\mu $ on 
$2^{X}$, a \textit{best transitive extension of }$\succsim $\textit{\
relative to} $D^{\mu }$ is any preorder $\succsim ^{\ast }$ $\in $ Ext$%
(\succsim )$ such that%
\begin{equation*}
D^{\mu }(\succsim ,\succsim ^{\ast })=\min \{D^{\mu }(\succsim
,\trianglerighteq ):\text{\thinspace }\trianglerighteq \text{\thinspace }\in 
\text{\thinspace Ext}(\succsim )\}\text{.}
\end{equation*}%
Fortunately, such extensions have a nice characterization.

\bigskip

\noindent \textbf{Theorem 5.1.} \textit{Let }$\mu $ \textit{be a positive
measure }$\mu $ \textit{on }$X$.\textit{\ Then, the unique best transitive
extension of any antisymmetric }$\succsim $ $\in \mathbb{A}(X)\ $\textit{%
with respect to }$D^{\mu }$\textit{\ is the transitive closure of} $\succsim 
$.

\bigskip

Before we prove this theorem, we present a simple example that shows that
the transitive closure of an antisymmetric acyclic order on $X$ need not be
a best approximation among all preorders on $X.$ This witnesses the
nontriviality of the general approximation we outlined above.

\bigskip

\noindent \textit{Example 5.1.} Let $X:=\{x_{1},x_{2},x_{3},x_{4}\},$ and
let $\succsim $ be the antisymmetric acyclic order on $X$ whose asymmetric
part is given as $x_{i}\succ x_{i+1}$ for $i=1,2,3.$ (The transitive closure 
$\succsim $ is the linear order on $X$ that ranks $x_{1}$ the highest, $%
x_{2} $ the second highest, so on.) Consider the reflexive binary relation $%
\trianglerighteq $ on $X$ whose asymmetric part is given as $%
x_{1}\vartriangleright x_{2}$ and $x_{3}\vartriangleright x_{4}$. Clearly, $%
\trianglerighteq $ is a partial order on $X,$ although it is not an
extension of $\succsim $. Moreover, 
\begin{equation*}
D(\succsim ,\trianglerighteq )=4<5=D(\succsim ,\text{tran}(\succsim )),
\end{equation*}%
so tran$(\succsim )$ is not a best approximation to $\succsim $ \textit{in }$%
\mathbb{P}(X).$ $\Vert $

\bigskip

We now turn to the proof of Theorem 5.1. Let us first observe that for any
antisymmetric $\succsim $ $\in \mathbb{A}(X),$ tran$(\succsim )$ is a
partial order on $X$ that extends $\succsim $ .\footnote{%
We use the antisymmetry postulate in Theorem 5.1 only to ensure that tran$%
(\succsim )$ is an antisymmetric extension of $\succsim $. As such, Theorem
5.1 applies to all non-antisymmetric $\succsim $ $\in \mathbb{A}(X)$ such
that tran$(\succsim )$ $\in $ Ext$(\succsim )$.
\par
Incidentally, note that tran$(\succsim )$ need not be an extension of a
reflexive relation $\succsim $ on $X$ that is either cyclic or not
antisymmetric. To illustrate, let $X:=\{a,b,c\}.$ If $\succsim $ equals $%
\Delta _{X}\sqcup \{(a,b),(b,c),(c,a)\},$ then $\succsim $ is a reflexive
and antisymmetric, but not acyclic, binary relation on $X,$ and tran$%
(\succsim )=X\times X$ which is not an extension of $\succsim $. On the
other hand, if $\succsim $ equals $(X\times X)\backslash \{(c,b)\},$ then $%
\succsim $ $\in \mathbb{A}(X)$ (but $\succsim $ is not antisymmetric) and
again tran$(\succsim )=X\times X$ which is not an extension of $\succsim $.}

\bigskip

\noindent \textbf{Lemma 5.2.} \textit{For any antisymmetric }$\succsim $%
\textit{\ }$\in \mathbb{A(}X)$\textit{, }tran$(\succsim )$ \textit{is a
partial order on }$X.$\textit{\ Moreover, for any} $\trianglerighteq $ $\in 
\mathbb{\ }$Ext$(\succsim ),$ \textit{we have }%
\begin{equation}
\succ \text{ }\subseteq \text{ tran}(\succsim )^{>}\subseteq \text{ }%
\vartriangleright \text{.}  \label{tra}
\end{equation}

\begin{proof}
Suppose $x$ tran$(\succsim )$ $y$ tran$(\succsim )$ $x$ for some distinct $%
x,y\in X.$ Then, there exist finitely many (pairwise distinct) $%
z_{0},...,z_{k},w_{0},...,w_{l}\in X$ such that $x=z_{0}\succsim
z_{1}\succsim \cdot \cdot \cdot \succsim z_{k}=y=w_{0}\succsim w_{1}\succsim
\cdot \cdot \cdot \succsim w_{l}=x.$ Since $\succsim $ is antisymmetric,
each $\succsim $ must hold strictly here, so we contradict acyclicity of $%
\succsim $. We thus conclude that tran$(\succsim )$ is antisymmetric, and
hence, a partial order on $X.\footnote{%
We give this argument here only for the sake of completeness. It is
well-known that an antisymmetric binary relation on a finite set is acyclic
if and only if its transitive closure is a partial order; see, for instance, 
\cite[Theorem 2.23]{Caspard}$.$}$

To prove (\ref{tra}), note that, by definition, $\succ $ $\subseteq $ tran$%
(\succsim )$. To derive a contradiction, suppose there exist $x,y\in X$ such
that $x\succ y$ but $y$ tran$(\succsim )$ $x$. Then, there exist an integer $%
k\geq 2$ and (pairwise distinct) $z_{0},...,z_{k}\in X$ with $%
y=z_{0}\succsim z_{1}\succsim \cdot \cdot \cdot \succsim z_{k}=x.$ Since $%
\succsim $ is antisymmetric, each $\succsim $ holds strictly, so we find $%
y\succ z_{1}\succ \cdot \cdot \cdot \succ x\succ y,$ contradicting the
acyclicity of $\succsim $. This proves the first containment in (\ref{tra}).
Next, suppose $x$ tran$(\succsim )^{>}$ $y.$ Then, again by antisymmetry of $%
\succsim ,$ there exist finitely many $z_{0},...,z_{k}\in X$ with $%
x=z_{0}\succ z_{1}\succ \cdot \cdot \cdot \succ z_{k}=y$. As $%
\trianglerighteq $ extends $\succsim $, we then have $x\vartriangleright
z_{1}\vartriangleright \cdot \cdot \cdot \vartriangleright $ $y,$ so, since $%
\trianglerighteq $ is transitive, we find $x\vartriangleright y$. This
proves the second containment in (\ref{tra}).
\end{proof}

\noindent \textbf{Lemma 5.2.} \textit{Let }$\succsim $\textit{\ be an
antisymmetric acyclic order on }$X\in $ \textit{and} $\trianglerighteq $ $%
\in \mathbb{\ }$Ext$(\succsim ).$ \textit{Then, }tran$(\succsim )$\textit{\
is in-between }$\succsim $\textit{\ and} $\trianglerighteq $.

\begin{proof}
The proof is by induction on the cardinality of the set tran$(\succsim
)\backslash \succsim $, say, $m.$ Consider first the case $m=1.$ Then, tran$%
(\succsim )\backslash \succsim $ $=\{(a,b)\}$ for some $a,b\in X.$ In view
of Lemma 5.2, $b\succ a$ cannot hold, so we have $(a,b)\in $ Inc$(\succsim
). $ Moreover, $a$ and $b$ are distinct (because $\succsim $ is reflexive)
so we have $a$ tran$(\succsim )^{>}$ $b$ (because tran$(\succsim )$ is
antisymmetric by Lemma 5.2). Again by Lemma 5.2, therefore, $%
a\vartriangleright b.$ It follows that tran$(\succsim )=$ $\succsim \oplus
(a,b)$ and $\succsim $\thinspace $\rightarrow $\thinspace tran$(\succsim )$%
\thinspace $\twoheadrightarrow $\thinspace $\trianglerighteq $, which means
tran$(\succsim )$\textit{\ }is in-between $\succsim $ and $\trianglerighteq $%
.

Now assume that $\succsim $\thinspace $\rightarrow $\thinspace tran$%
(\succsim )$\thinspace $\twoheadrightarrow $\thinspace $\trianglerighteq $
holds for every antisymmetric $\succsim $\textit{\ }$\in \mathbb{A}(X)$ and $%
\trianglerighteq $ $\in \mathbb{\ }$Ext$(\succsim )$ such that tran$%
(\succsim )\backslash \succsim $ has $m\geq 1$ elements. To complete the
induction, suppose $\succsim $ is an antisymmetric acyclic order on $X$ with 
$\left\vert \text{tran}(\succsim )\backslash \succsim \right\vert =m+1.$
Pick any $(a,b)$ in tran$(\succsim )\backslash \succsim $. By the same
argument we made in the previous paragraph, we must have $(a,b)\in $ Inc$%
(\succsim )$ and $a\vartriangleright b.$ Moreover, acyclicity of $\succsim $
entails that of $\succsim _{0}$ $:=$ $\succsim \sqcup \{(a,b)\}.$ (For,
otherwise, there exist finitely many $z_{1},...,z_{k}\in X$ with $z_{1}\succ
\cdot \cdot \cdot \succ z_{k}\succ z_{1}$. Since $\succsim $ is acyclic, $%
(z_{i},z_{i+1})=(a,b)$ for some $i=1,...,k-1,$ and we can take $i=1$,
relabelling if necessary. But since $a$ tran$(\succsim )$ $b,$ there also
exist finitely many $w_{0},...,w_{l}\in X$ with $a=w_{0}\succ \cdot \cdot
\cdot \succ w_{l}=b$. Consequently, $b=z_{2}\succ \cdot \cdot \cdot \succ
z_{k}\succ z_{1}=a=w_{1}\succ \cdot \cdot \cdot \succ w_{l}=b,$
contradicting the acyclicity of $\succsim $.) Thus, tran$(\succsim _{0})=$ $%
\succsim \oplus (a,b)$ and $\succsim $\thinspace $\rightarrow $\thinspace $%
\succsim _{0}$\thinspace $\twoheadrightarrow $\thinspace $\trianglerighteq $%
. Now notice that tran$(\succsim _{0})\backslash \succsim _{0}$ has $m$ many
elements, so by the induction hypothesis, $\succsim _{0}$\thinspace $%
\rightarrow $\thinspace tran$(\succsim _{0})$\thinspace $\twoheadrightarrow $%
\thinspace $\trianglerighteq $. It follows that $\succsim $\thinspace $%
\rightarrow $\thinspace tran$(\succsim _{0})$\thinspace $\twoheadrightarrow $%
\thinspace $\trianglerighteq $. Since tran$(\succsim )=$ tran$(\succsim
_{0}),$ we are done.
\end{proof}

Now let $\mu $ be a positive measure $\mu $ on $X$, take any antisymmetric $%
\succsim $ $\in \mathbb{A}(X),$ and let $\trianglerighteq $ $\in \mathbb{\ }$%
Ext$(\succsim )$. Then, by Lemma 5.2, tran$(\succsim )$ is in-between $%
\succsim $ and $\trianglerighteq $. As $D^{\mu }$ satisfies Axiom 1, we thus
get%
\begin{equation*}
D^{\mu }(\succsim ,\trianglerighteq )=D^{\mu }(\succsim ,\text{tran}%
(\succsim ))+D^{\mu }(\text{tran}(\succsim ),\trianglerighteq )\geq D^{\mu
}(\succsim ,\text{tran}(\succsim ))\text{.}
\end{equation*}%
This completes the proof of Theorem 5.1.

\section{Embedding $(\mathbb{A}(X),D)$ in a Euclidean Space}

As $(\mathbb{A}(X),D^{\mu })$ is a finite metric space (for any positive
measure on $\mu $ on $2^{X}$), another natural query to consider here is if
we can embed this space isometrically in a Euclidean space. There are of
course well-known characterizations of finite metric spaces that are
isometrically embeddable in a Euclidean space; see, for instance, \cite%
{Bowers,Morgan,Sch}. But even without resorting to such theorems, we can
show easily that $(\mathbb{A}(X),D^{\mu })$ is not isometrically embeddable
in a Euclidean space. In fact, endowing the set $\mathbb{P}^{\ast }(X)$ of
all partial orders on $X$ with $D^{\mu }$ yields a metric space which is not
isometric to any subset of a Euclidean space (unless $X$ contains only two
elements). More generally:

\bigskip

\noindent \textbf{Proposition 6.1.} \textit{Let }$X$ \textit{be a finite set
with }$\left\vert X\right\vert \geq 3$\textit{\ and }$\mu $ \textit{a
positive measure on }$2^{X}$\textit{. Then, }$(\mathbb{P}^{\ast }(X),D^{\mu
})$ \textit{cannot be isometrically embedded in any strictly convex }(%
\textit{real})\textit{\ normed linear space.}

\begin{proof}
Take any distinct $a,a^{\prime },b\in X$ and define $\succsim _{0}$ $%
:=\Delta _{X}\sqcup \{(a,b)\},$ $\succsim _{1}$ $:=\Delta _{X}\sqcup
\{(a^{\prime },b)\}$ and $\trianglerighteq $ $:=\Delta _{X}\sqcup
\{(a,b),(a^{\prime },b)\}.$ These are partial orders on $X,$ and it is plain
that $\succsim $\thinspace $\rightarrow $\thinspace $\succsim _{0}$%
\thinspace $\twoheadrightarrow $\thinspace $\trianglerighteq $ and $\succsim 
$\thinspace $\rightarrow $\thinspace $\succsim _{1}$\thinspace $%
\twoheadrightarrow $\thinspace $\trianglerighteq $. Now, to derive a
contradiction, let us suppose that there is an isometric embedding $\varphi :%
\mathbb{P}^{\ast }(X)\rightarrow E$ for some strictly convex normed linear
space $E$ (whose norm is denoted as $\left\Vert \cdot \right\Vert $). By
translation, we may assume $\varphi (\Delta _{X})=\mathbf{0}_{E},$ where $%
\mathbf{0}_{E}$ is the origin of $E$. As $D^{\mu }$ satisfies Axiom 1, and $%
\varphi $ is an isometry, we have%
\begin{eqnarray*}
\left\Vert \varphi (\succsim _{0})\right\Vert +\left\Vert \varphi
(\trianglerighteq )-\varphi (\succsim _{0})\right\Vert &=&D^{\mu }(\Delta
_{X},\succsim _{0})+D^{\mu }(\succsim _{0},\trianglerighteq ) \\
&=&D^{\mu }(\Delta _{X},\trianglerighteq ) \\
&=&\left\Vert \varphi (\trianglerighteq )\right\Vert .
\end{eqnarray*}%
Since $E$ is strictly convex, therefore, there exists a positive real number 
$\lambda _{0}$ such that $\varphi (\succsim _{0})=\frac{\lambda _{0}}{%
1+\lambda _{0}}\varphi (\trianglerighteq ).$ Precisely the same reasoning
yields also that $\varphi (\succsim _{1})=\frac{\lambda _{1}}{1+\lambda _{1}}%
\varphi (\trianglerighteq )$ for some $\lambda _{1}>0.$ Since 
\begin{equation*}
\left\Vert \varphi (\succsim _{0})\right\Vert =D^{\mu }(\Delta _{X},\succsim
_{0})=2^{\left\vert X\right\vert -1}\mu (\{b\})=D^{\mu }(\Delta
_{X},\succsim _{1})=\left\Vert \varphi (\succsim _{1})\right\Vert ,
\end{equation*}%
we must have $\lambda _{0}=\lambda _{1},$ and it follows that $\varphi
(\succsim _{0})=\varphi (\succsim _{1}),$ whence $\succsim _{0}$ $=$ $%
\succsim _{1}$, a contradiction.
\end{proof}

It follows from this result that any sort of embedding of $(\mathbb{P}^{\ast
}(X),D^{\mu }),$ and hence of $(\mathbb{P}(X),D^{\mu })$ or $(\mathbb{A}%
(X),D^{\mu }),$ in a Euclidean space must involve some distortion.\footnote{%
The \textit{distortion} of a map $\varphi $ from a metric space $(Z,d)$ into
a normed linear space $E$ is defined as the product of $\sup_{w,z\in Z}\frac{%
\left\Vert \varphi (w)-\varphi (z)\right\Vert }{d(w,z)}$ and $\sup_{w,z\in Z}%
\frac{d(w,z)}{\left\Vert \varphi (w)-\varphi (z)\right\Vert }.$} For
instance, by a famous theorem of Bourgain \cite{Bourgain}, $(\mathbb{A}%
(X),D^{\mu })$ can be embedded in the Hilbert space $\ell _{2}$ with
distortion at most $O(\log \left\vert \mathbb{A}(X)\right\vert ).$ We shall
not pursue this problem in this paper any further, however.

\section{Future Research}

In conclusion, we would like to point out some directions for future
research. First, there are some natural best approximation problems that one
should attack. A really interesting one, for instance, concerns finding the
nearest total preorder on $X$ to any given preorder $\succsim $ on $X$ in
terms of the metric $D.$ This sort of a study would aim at characterizing
such best complete approximations algebraically as well as algorithmically.
This may be particularly useful when the incompleteness of a preference
relation arises due to \textquotedblleft missing data.\textquotedblright\
Moreover, it would allow approximating various decision problems and games
with incomplete preferences by more standard models. In addition, it would
furnish a natural index of incompleteness, namely, the minimum $D$ distance
between $\succsim $ and its projection onto the set of all complete
preorders on $X.$

Second, one may take up the problem of deducing consensus preferences from a
given family of preferences, say, by minimizing the sum of $D$ distances
from that family. These sorts of problems are NP-hard, and studied
extensively in the operations research literature in terms of the
Kemeny-Snell-Bogart metric. It should be very interesting to find out the
consequences of replacing $d_{\text{KSB}}$ with $D$ in those studies.

Finally, we should note that the majority of economic models presume
infinite alternative spaces, and indeed, the most well-known models of
individual decision theory, such as the expected utility model under risk
and uncertainty, the model of Knightian uncertainty, time discounting
models, menu preferences, etc., work with preferences that are defined on an
infinite alternative space. By contrast, our work in this paper depends very
much on the finiteness of $X,$ and while it is readily applicable to
experiments, individual choice theory, voting, etc., it does not play well
within these settings. One may, of course, always extend the top-difference
semimetric $D$ to the case of an arbitrary $X$ by means of the formula%
\begin{equation*}
D(\succsim ,\trianglerighteq )=\sup \sum_{S\subseteq X}\left\vert
M(S,\succsim )\triangle M(S,\trianglerighteq )\right\vert ,
\end{equation*}%
where $\sup $ is taken over all finite subsets of $X,$ but this seems like a
rather coarse approach. (It would not, for instance, distinguish any two
quasi-linear preferences on $\mathbb{R}^{2}$.) Extending the approach
developed here to the context of infinite alternative spaces remains as a
major problem for future research.

\bigskip

\noindent \textsc{Acknowledgement }We thank professors Chris Chambers and
Federico Echenique for their insightful comments in the development stage of
this work.

\end{document}